\documentclass[letter,nofootinbib, superscriptaddress,twocolumn]{revtex4-1}
\usepackage{amsmath}
\usepackage{amssymb}
\usepackage{braket}
\usepackage{amssymb}
\usepackage{color}
\usepackage{xcolor}
\usepackage{graphicx}
\usepackage{latexsym}
\usepackage{listings}
\usepackage{mathrsfs}

\usepackage{subfigure}
\usepackage{verbatim}
\usepackage{tabularx}
\usepackage{ragged2e}
\usepackage{multirow}
\usepackage{amsbsy}
\usepackage{siunitx}
\usepackage{numprint}
\usepackage{makecell}
\usepackage[utf8]{inputenc}

\usepackage[colorlinks,citecolor=blue]{hyperref}

\newcommand{\paren}[1]{\left( #1 \right)}

\newcommand{\inner}[1]{\paren{#1}}
\newcommand{\norm}[1]{\left| #1 \right|_2}

\renewcommand{\braket}[1]{\left\langle #1 \right\rangle}

\newcommand{\V}{\mathcal{V}}

\newcommand{\M}{\mathcal{M}}
\newcommand{\DN}{\mathcal{N}}

\newcommand{\Ord}[1]{\mathcal{O}\paren{#1}}

\newcommand{\revone}[1]{#1}%{\color{magenta}#1}}
\newcommand{\revtwo}[1]{#1}%{{\color{orange}#1}}

\npdecimalsign{.}
\npthousandsep{,}

\newcolumntype{Y}{>{\RaggedRight\arraybackslash}X}

\sloppypar

\begin{document}

\title{Spectral Methods in the Presence of Discontinuities}

\author{Joanna Piotrowska}
\email{jmp218@cam.ac.uk}
\affiliation{Cavendish Laboratory, University of Cambridge, Cambridge, UK}
\affiliation{Perimeter Institute for Theoretical Physics, Waterloo,
  ON, Canada}

\author{Jonah M. Miller}
\email{jonahm@lanl.gov}
\affiliation{Computational Physics and Methods, Los Alamos National Laboratory,
  Los Alamos, NM, USA}
\affiliation{Center for Nonlinear Studies, Los Alamos National Laboratory,
  Los Alamos, NM, USA}
\affiliation{Center for Theoretical Astrophysics, Los Alamos National Laboratory,
  Los Alamos, NM, USA}

\author{Erik Schnetter}
\email{eschnetter@perimeterinstitute.ca}
\affiliation{Perimeter Institute for Theoretical Physics, Waterloo,
  ON, Canada}
\affiliation{Department of Physics and Astronomy, University of
  Waterloo, Waterloo, ON, Canada}
\affiliation{Center for Computation \& Technology, Louisiana State
  University, Baton Rouge, LA, USA}

\date{2018-12-11}

\begin{abstract}

  Spectral methods provide an elegant and efficient way of numerically
  solving differential equations of all kinds. For smooth problems,
  truncation error for spectral methods vanishes exponentially in the
  infinity norm and $L_2$-norm. However, for non-smooth problems,
  convergence is significantly worse---the $L_2$-norm of the error for
  a discontinuous problem will converge at a sub-linear rate and the
  infinity norm will not converge at all. We explore and improve upon
  a post-processing technique---optimally convergent mollifiers---to
  recover exponential convergence from a poorly-converging spectral
  reconstruction of non-smooth data. This is an important first step
  towards using these techniques for simulations of realistic systems.

\end{abstract}

\maketitle

\begin{figure*}[t]
  \centering
  \includegraphics[width=\textwidth]{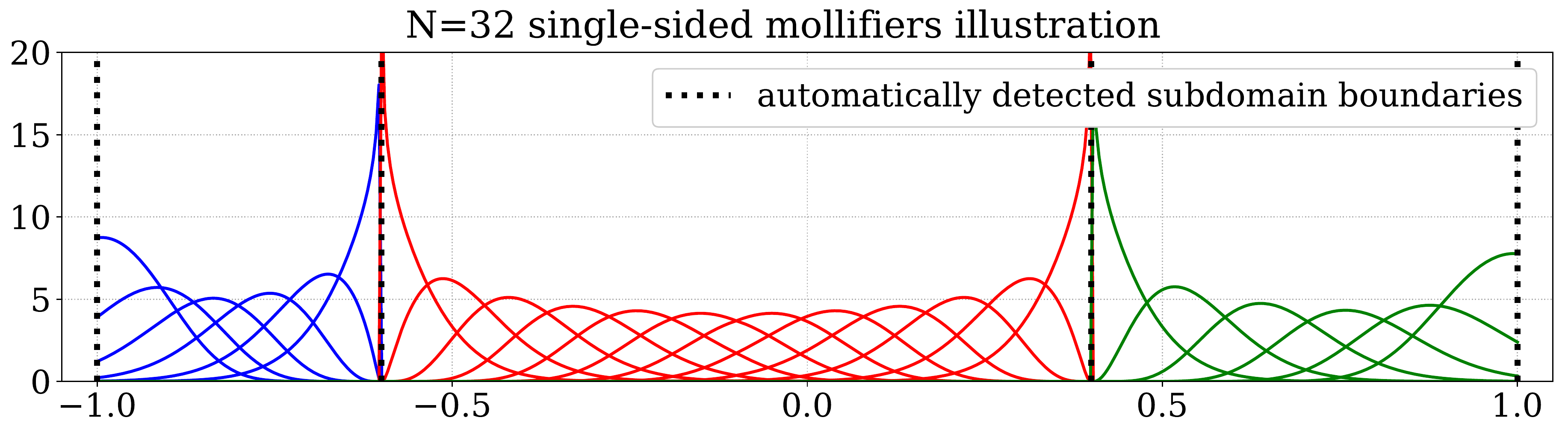}
  \caption{Illustration of the shapes of adaptive single-sided
    mollifiers both near and far from a boundary for piecewise
    mollification.  Figure presents a selected sample of 24 mollifiers
    constructed for a piecewise flat function on the domain $[-1,1]$
    with discontinuities at $x=-0.6$ and $x=0.4$. The mollifier's
    support varies in space, becoming most compact close to the
    function discontinuity. In order to avoid combining information
    from two regions of smoothness across the discontinuity, the
    mollifying kernel is non-zero only on one side of the
    discontinuity.}
  \label{fig:mollifiers}
\end{figure*}

Spectral methods are a means of approximating arbitrary, integrable
functions in the discrete language of a computer. They provide an
alternative to the widely applied finite difference and finite
volume-type methods.\footnote{One can think of a finite element method
  as a multi-domain spectral method.} Because (for smooth functions
and in the appropriate norm) they enjoy exponential convergence onto
true solution across the entire domain, they are a powerful tool for
solving partial differential equations of all kinds.\footnote{Here we
  discuss the application of spectral methods to hyperbolic problems
  via the method of lines. However, we believe our analysis has broad
  applicability to other problem types.} For non-smooth problems, such
as the modeling of fluid shocks, spectral methods typically perform
poorly. They suffer spurious oscillations near discontinuities and, in
the worst case, converge slower than linearly. This is a manifestation
of the famous \textit{Gibbs phenomenon}
\cite{wilbraham1848certain,gibbs1898fourier,MichelsonGibbsReply,boyd2013chebyshev}.\footnote{For
  a historical perspective on the Gibbs phenomenon, see Hewitt and
  Hewitt \cite{Hewitt1979}.}

Because of the extreme efficiency of spectral methods (and indeed,
their great success in many domains such as relativistic
astrophysics---see, e.g., \cite{BoyleBBH,ScheelWaveforms,SPeC}), it is
desirable to extend these techniques to non-smooth problems. Many
authors over the years have attempted to evade the Gibbs phenomenon in
the non-smooth case and recover the exponential convergence that makes
spectral methods so powerful. Approaches include (but are not limited
to) filtering \cite{Vandeven1991a}, artificial viscosity
\cite{Tadmor1990}, reprojection \cite{Gottlieb1992a,Gottlieb2011a},
and mollification \cite{gottlieb1985recovering}.\footnote{For a review of
  mollifiers, see \cite{Tadmor2007}.}

One especially intriguing way of dealing with discontinuities is to
\textit{model} and then \textit{subtract} the error introduced by the
Gibbs phenomenon. This idea has been explored in various contexts,
perhaps starting in 1906 with Krylov \cite{krylov2006approximate},
to be later revisited by Lanczos in 1966 \cite{lanczos2016discourse}
and again by Eckhoff  in 1994 \cite{ECKHOFF1994103}. A particularly
flexible modern iteration of this approach, first introduced by Lipman
and Levin, uses moving least squares to fit the error introduced by
the discontinuity \cite{LipmanLevinLeastSquares,AMIR201831}.

\revtwo{Alternatively, if the positions of discontinuities are known ahead of time, a
\textit{global} spectral method may be replaced by a
\textit{multi-domain} spectral method, where discontinuities lie on
domain boundaries
\cite{canuto1982approximation,funaro1988new,funaro1991convergence}. Although
it is not usually presented in this way, the most popular of these
approaches is no-doubt the discontinuous Galerkin method
\cite{hesthaven2007nodal}.}\footnote{For a discussion on how the
  discontinuous Galerkin method fits into the more broad family of
  multi-domain spectral methods, see
  \cite{GottliebSpectralForHyperbolic} and references therein.}

% These
% methods are extraordinarily powerful, however since they require the
% location of discontinuities to be known, we will focus our attention
% elsewhere in this work.

Although spectral evasions of the Gibbs phenomenon have been studied
extensively for decades, they have barely made it into practical
applications. (With a few notable exceptions---see,
e.g., \cite{Gelb2000Enhanced,MeisterFilter}.)
%\TODO{Find more references. -JMM} 
In this work, we explore and improve upon the spectral mollifiers
developed by Tadmor, \revone{Tanner}, and collaborators
\cite{gottlieb1985recovering}, \cite{Tadmor2002}, \cite{Tanner2006}
with the intent of applying these techniques to relativistic
astrophysics. In particular, we develop the family of
\textit{one-sided} mollifiers, which incorporate the discontinuous
nature of the underlying function. \revone{A distinguishing feature of
mollifiers} is that there is a different mollifying function for every
point in the domain.  Figure~\ref{fig:mollifiers} illustrates this
fact with a representative set of these functions for a function with
discontinuities at $x=-0.6$
$x=0.4$. % \TODO{This is added so we reference figure 1 first. -JMM}

We combine these mollifiers with
the spectral edge detection developed by Gelb and collaborators
\cite{Gelb1999a,Gelb2001a,Gelb2008a}. We offer practical advice
on the subtleties of implementing these tools, discuss their
limitations and compare mollification to the Gegenbauer 
reconstruction.
%We compare mollification to the
%Gegenbauer reconstruction, and We offer practical advice on the
%subtleties of implementing these tools and we discuss their
%limitations. 
Although mollifiers and edge detectors are established technology,
they have rarely, if ever, been applied in real situations. We
therefore believe that the improvements and practical advice we offer
in this paper comprise an important first step towards applying these
powerful techniques in realistic, non-smooth situations.

In Section \ref{sec:intro}, we discuss spectral methods and their
limitations and introduce spectral edge detection and adaptive
mollifiers. In Section \ref{sec:decomposition}, we specialize our
discussion to the Chebyshev basis. In Section
\ref{sec:edge:detection}, we discuss edge detection in more detail and
present several numerical experiments that show its strengths and
weaknesses. In Section \ref{sec:mollifiers}, we discuss spectral
mollification and present several experiments. In Section
\ref{sec:hyperbolic:PDEs} we integrate edge detection and
mollification to solve the linear advection equation with
discontinuous initial data. We compare our results to the Gegenbauer
reconstruction and show that even in the case of simplest
implementation, mollification behaves in a more robust
fashion. Finally, in Section \ref{sec:conclusion}, we offer some
concluding thoughts and describe the path forward.

% When combined with the method of lines, spectral methods can solve a
% well-posed initial value problem by discretizing spatial dependence of
% a function into a sum of basis functions (often interpolating
% polynomials) and reducing the problem to a set of ordinary
% differential equations. These can be solved using standard numerical
% integration techniques. \TODO{Spectral methods are not limited to a
%   method-of-lines approach. You can discretize the \textit{time
%     direction} in the same way. (Spectral methods also have broad
%   applicability in boundary-value problems. Above is my best attempt
%   at dancing around these technicalities. -JMM}

%\vfill
%\pagebreak

\section{Spectral Methods}
\label{sec:intro}

\subsection{Pseudospectral projection}
\label{sec:intro:projection}

The principle idea behind spectral methods (and indeed most numerical
methods) is to represent a function of interest $u(x)$ as a linear
combination of a set of suitably chosen basis functions $\phi_n(x)$:
\begin{equation}
\label{eq:S_N}
u(x) \approx S_N[u](x) = \sum_{n=0}^{N} \hat{u}_n \phi_n(x),
\end{equation}
where $S_N[u]$ is the partial sum of $u$, the expansion coefficients
are given by
\begin{equation}
\label{eq:f_n}
\hat{u}_n=\frac{\inner{u_n,\phi_n}}{\inner{\phi_n, \phi_n}},
\end{equation}
and
\begin{equation}
  \label{eq:inner:product}
  \inner{a,b} = \int _\Omega a(x) b(x) w(x) dx
\end{equation}
is the inner product between $a$ and $b$ over some domain $\Omega$
with weight function $w(x)$.

For smooth functions, this representation suffers an error whose
infinity norm decays faster than any power of $N$ and which vanishes
in the limit of $N\rightarrow\infty$
\cite{GrandclementSpectral,GottliebSpectralForHyperbolic}.\footnote{The
  details of convergence can vary with spectral representation and
  with both the smoothness and analyticity of the function being
  approximated. For example, the spectral representation of an
  analytic function will converge more rapidly than a representation
  of a smooth one, although both will converge exponentially. For a
  formal overview of these issues, see \cite{boyd2013chebyshev}.}
Handling nonlinear differential equations in the basis-coefficient (or
\textit{modal}) representation requires regular integration over the
whole domain, which is unfavorable. Instead, one can choose a set of
\textit{collocation points}\footnote{The optimal choice of collocation
  points depends on the spectral basis chosen.} and use the Lagrange
polynomials
\begin{equation}
  \label{eq:lagrange}
  L_j(x) = \prod_{\substack{0<m<k\\m\neq j}}\frac{x-x_m}{x_j-x_m}
\end{equation}
as a second set of basis functions. In this case, integration can
be approximated via Gauss quadrature:
\begin{equation}
\label{gauss_quad}
\int_\Omega f(x)w(x)dx \approx \sum_{n=0}^{N} f(x_n) w_n ,
\end{equation}
where $x_n$ are chosen collocation points and $w_n$ are their
associated weights. In this approach, $S_N[u]$ \textit{interpolates}
$u$ between the collocation points $x_n$. This family of methods are
known as \textit{pseudospectral} or \textit{collocation-spectral}
methods.

One can easily transform between the collocation representation
$u_j = u(x_j)$ of a function (which is just the Lagrange basis in
disguise) and the spectral coefficients $\hat{u}_i$ of its
interpolating partial sum via the matrix operation
\begin{equation}
\label{V}
\hat{u}_i = \V _{ij}u_j, \qquad \V_{ij} \equiv
    \frac{\phi_i(x_j)w_j}{(\phi_i,\phi_i)},
\end{equation} 
where $\V_{ij}$ is the inverse of the \textit{Vandermonde
  matrix}. Similarly, one can approximately map the partial sum
$S_N[u]$ of $u$ to the partial sum $S_N[\partial_x u]$ of its
derivative via the \textit{modal} differentiation matrix
\begin{equation}
  \label{eq:D:modal}
  \M_{ij} = \frac{\inner{\partial_x \phi_i,\phi_j}}
    {\inner{\phi_j,\phi_j}}.
\end{equation}
Moreover, one can approximately map the restrictions $u_i$ to
derivative of $u$ restricted to $x_i$ via the \textit{nodal}
differentiation matrix
\begin{equation}
  \label{eq:D:nodal}
  \DN = \V^{-1} \M \V,
\end{equation}
which is nothing more than the modal differentiation matrix
appropriately transformed into the collocation representation
\cite{GrandclementSpectral,GottliebSpectralForHyperbolic}. This
results in the following convenient approximations:
\begin{equation}
\label{Diff_matrix}
\frac{d}{dx}u(x_j) \approx \DN_{ji}u(x_i), \qquad 
    \frac{d}{dx}\hat{u}_i \approx \M_{ij}\hat{u}_j
\end{equation}

\subsection{Influence of discontinuities}
\label{sec:intro:gibbs}

\begin{figure*}
  \centering
  \subfigure[]{
  \label{fig:gibbs:nogibbs}
  \includegraphics[width=.48\textwidth]{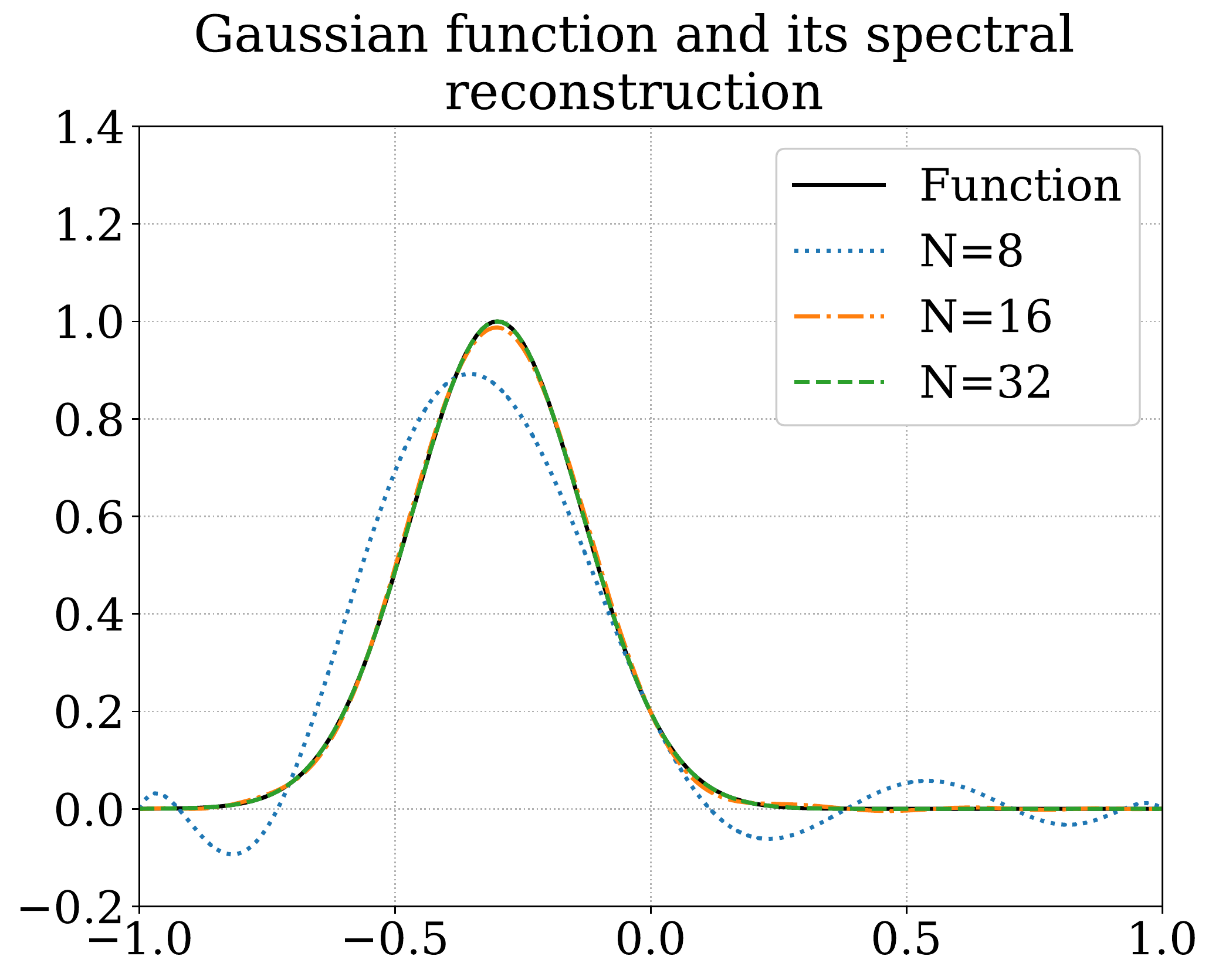}
  }
  \subfigure[]{
  \label{fig:gibbs:gibbs}   
   \includegraphics[width=.48\textwidth]{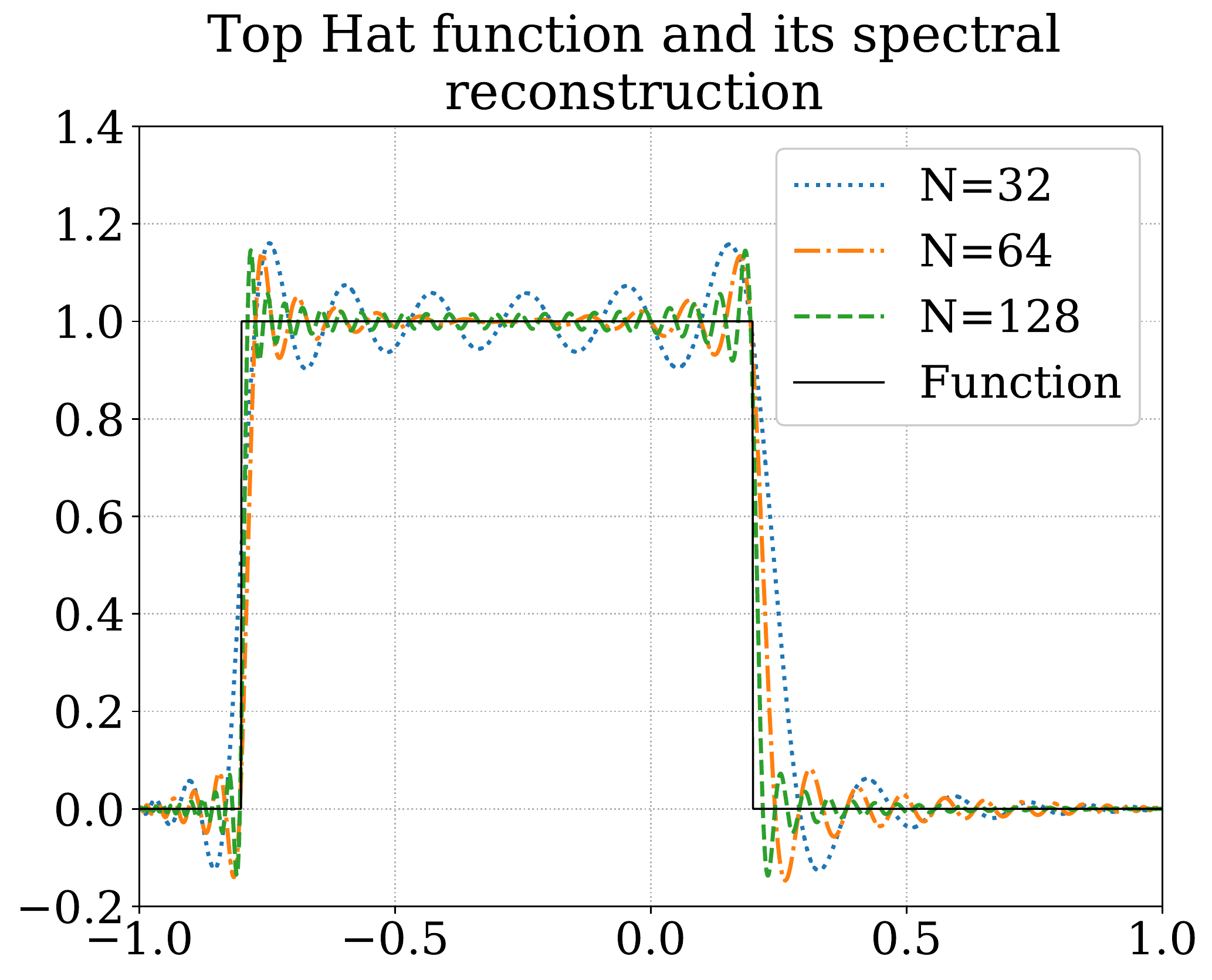}
   }
  \caption{Illustration of spectral reconstruction and the influence 
    of discontinuities: \subref{fig:gibbs:nogibbs} exponentially 
    convergent expansion of a~smooth Gaussian. At expansion order 
    $N=16$ the reconstruction already resembles the true solution. 
    \subref{fig:gibbs:gibbs} Gibbs phenomenon on a simple example of 
    a Top Hat function. High frequency oscillations near 
    discontinuities do not vanish with increasing number of modes 
    $N$, as they retain constant amplitude and get `pushed closer' 
    towards the discontinuity.}
  \label{fig:gibbs}
\end{figure*}

Unfortunately, the spectral expansion is only so extraordinarily
accurate for smooth functions. For Chebyshev polynomials, the $L_2$-norm 
of the difference between an arbitrary $C^m$ function $u$ and its 
partial sum $S_N[u]$ is approximately given by
\begin{equation}
  \label{eq:error:bound:chebyshev}
  \norm{u - S_N[u]} \leq \frac{\alpha}{N^{m}}\sum_{k=0}^m 
    \norm{f^{(k)}}
\end{equation}
for some positive constant $\alpha$
\cite{GrandclementSpectral,GottliebSpectralForHyperbolic}. And one can
draw similar bounds for other sets of orthogonal basis functions.
This means that, in presence of discontinuities, convergence onto the
function of interest is non-uniform and pointwise slow---the error
decays as $\Ord{N^{-1}}$ away from the discontinuity and the $L_2$
norm exhibits at best $\Ord{N^{-\frac{1}{2}}}$ convergence
\cite{thomson2008elementary}. Worse, the infinity norm of the error
does not converge at all! It is the local behaviour in the vicinity of
the function's discontinuity which spoils the global approximation
within the domain---this is the famous \textit{Gibbs phenomenon}
\cite{wilbraham1848certain,gibbs1898fourier,MichelsonGibbsReply,boyd2013chebyshev}.

Close to the discontinuity, the reconstruction is polluted by
high-frequency oscillations of non-decreasing amplitude, irrespective
of the $N$ number of modes used, as shown in
Figure~\ref{fig:gibbs}. This behaviour prevents the recovery of true
solution, making the method unsuitable for discontinuous problems of
interest, such as fluid shocks or stellar surfaces.

\subsection{Recovering the underlying solution free from Gibbs 
phenomenon}
\label{sec:intro:mollification}
Before we introduce methods designed for eliminating the Gibbs 
phenomenon, we mention an important property of the Chebyshev
polynomials:
\begin{equation}
\label{eq:chebcos}
T_n(\cos\theta)=\cos(n\theta)
\end{equation}
which implies that for any function, $v:[-1,+1]\rightarrow\mathbb{R}$
a~change of variables $u(\theta) = v(\cos\theta)$ implies that
\begin{equation}
\label{eq:chebfour}
\int_{-1}^{+1}v(x)T_k(x){(1-x^2)}^{-\frac{1}{2}}dx=
\frac{2}{\pi}\int_0^{2\pi}u(\theta)\cos(k\theta) d\theta.
\end{equation}
Equations \eqref{eq:chebcos} and \eqref{eq:chebfour} tell us that
Chebyshev expansion is equivalent to a~Fourier cosine series with a 
change of variables.  This, in turn, allows us to study Fourier series 
and apply what we learn to their Chebyshev equivalents. In this 
subsection we shall refer to the analysis of Fourier expansion, in 
order to remain consistent with the literature.

\subsubsection{Filtering}
\label{sec:intro:mollification:filters}

Vandeven \cite{Vandeven1991a} suggests that convergence of 
a~Fourier series $S_N f(x)$ of a discontinuous function $f(x)$:
\begin{equation}
    S_N f(x) = \sum_{k=-N}^N \hat{f}(k)e^{ikx}
\end{equation}
can be accelerated by introducing spectral 
filters $\sigma(\frac{k}{N})$, which are smooth functions of compact 
support $[-1,1]$, characterised by $\sigma(0)=1$. Upon application
of $\sigma(\frac{k}{N})$, one arrives at filtered spectral expansion
of function $f(x)$, $\mathcal{S}^\sigma_N(x)$ defined below:
\begin{equation}
\label{eq:van:filters}
    \mathcal{S}^\sigma_N(x)= \sum_{k=-N}^N \hat{f}(k) \sigma\Big(\frac{k}{N}\Big)e^{ikx}
\end{equation}
which converges to $f(x)$ faster than the Gibbs-affected $S_N f(x)$. 
This is achieved by decreasing the importance of higher order terms 
in the spectral expansion, which smooths out Gibbs oscillations, 
while preserving the low-frequency bulk approximation to $f(x)$. 
However, this operation is equivalent to real space convolution 
with a~Fourier transform of the $\sigma(\frac{k}{N})$ filter, 
which acts like a~smoothing kernel and results in global smoothing 
of the function and subsequent smoothing out of the discontinuity.

One way to avoid this \textit{global} convolution and thus smoothing
of the discontinuity is to change the strength of the filter $\sigma$
as a function of space so that
\begin{equation}
  \label{eq:spacewise:filter}
  \sigma(k/N) \to \sigma(x,k/N),
\end{equation}
as in (for example) \cite{Boyd89Filter}. We take a related approach,
\textit{mollification}, which is described below.

\subsubsection{Mollification}
\label{sec:intro:mollification:mollifiers}

In \cite{gottlieb1985recovering}, Gottlieb and Tadmor use the idea of 
spectral filters to
introduce \textit{mollification}, which is a physical space equivalent
of \textit{filtering} carried out in Fourier space:
\begin{equation}
\label{eq:molli:filtering}
\Phi_{p,\delta}*(S_Nf)(x) \longleftrightarrow 
\sum_{|k| \leq N} \varphi_{p,\delta}\Big(\frac{k}{N}\Big)
    \hat{f}(k)e^{ikx},
\end{equation}
where $\varphi_{p,\delta}$ is equivalent to the $\sigma$ in equation
\eqref{eq:van:filters} above. They define mollifiers $\Phi_{p,\delta}$
as unit mass, non-negative kernels of compact support. These are then
adapted at every point within the domain, such that the essential
support---i.e., the majority of the nonzero portion---of
$\Phi_{p,\delta}*(S_Nf)(x)$ does not cross the discontinuities in $f$.
To avoid the discontinuities, the method must identify where they
are. Therefore, recovery of the Gibbs-free solution requires two basic
steps: \textit{edge detection}, followed by real-space
\textit{mollification}.

\revone{In \cite{Tanner2006}}, Tanner defines ``optimal'' mollifiers
by starting from a family of filters $\varphi(\xi)$ characterised by
optimal exponential decay in both physical and Fourier
space:\footnote{\revone{We follow the treatment summarized by Tadmor
    in the review article \cite{Tadmor2007}.}}
\begin{equation}
    \varphi(\xi):=\exp\bigg(-\frac{\xi^2}{2}\bigg) \times \Bigg[\sum_{j=0}^p 
    \frac{1}{2^j j!}\xi^{2j}\Bigg].
\end{equation}

% \pagebreak

He then takes their inverse Fourier transform
\begin{equation}
   \label{eq:mollification:mollifier:inverse:transform}
   \Phi(y):= \int_{I\!R}\exp\bigg(2 \pi i y \xi\bigg) \varphi (\xi) d\xi
\end{equation}
to define a global
mollifier $\Phi_p(y)$:
\begin{equation}
   \label{eq:mollification:mollifier:undilated}
      \Phi_p(y):=\exp\bigg(-\frac{y^2}{2}\bigg) \times \Bigg[\sum_{j=0}^p 
    \frac{(-1)^j}{4^j j!}H_{2j}\bigg(\frac{y}{\sqrt{2}}\bigg)\Bigg],
    \vspace{0.1cm}
\end{equation}
where $H_{\alpha}$ is the Hermite polynomial of order $\alpha$
and the number $p$ of frequencies used in the inverse Fourier
transform is decided by compactness concerns as described below.

%\vfill
%\newpage
%\vspace{2cm}

Finally, to ensure appropriate essential support, $\Phi(x)$ 
is dilated as
\begin{equation}
    \label{eq:dilation}
\Phi_\delta(x,y):= \frac{1}{\delta(x)}\Phi_p\bigg(\frac{y}{\delta(x)}\bigg),
\end{equation}
where $\delta(x)$ is given by
\begin{equation}
  \label{eq:def:delta}
    \delta(x):= \sqrt{\theta d(x) N}
\end{equation}
for a free parameter $\theta$.\footnote{In a numerical investigation
  performed in \revone{\cite{Tanner2006}}, Tanner uses the value of
  $\theta \sim 1/4$, which we implement in our tests as
  well. \revone{This is not necessarily the optimal choice.}} The
parameter $d(x)$, on the other hand, introduces information about the
discontinuities within the domain.
% \begin{equation}
%     d(x):=\mathrm{dist}\Big\{x,\big\{c_1,\ldots,c_J\big\}\Big\},
% \end{equation}
% where $c_1,\ldots, c_J$ are detected locations of the discontinuities
% and 
The interval $(x-d(x), x+d(x))$ is the largest interval of smoothness enclosing
and centered on
$x$.\footnote{A larger region of continuity could be constructed which
  encloses $x$ but is not centered on it. This could, in principle,
  improve the effectiveness of the mollifier and is worthy of future study.} The number of terms $p$ in the sum in equation
\eqref{eq:mollification:mollifier:undilated} is then given by:
\begin{equation}
  \label{eq:def:fourier:p}
    p=p_N:=\theta^2d(x)N.
\end{equation}

Using this recipe, Tanner arrives at the exponentially accurate
mollifier $\Phi^{d}_{p_N}(x,y)$  given by:
\begin{widetext}
\begin{equation}
\label{eq:the:mollifier}
    \Phi^{d}_{p_N}(x,y)=\frac{1}{\delta(x)}\exp\Big(-\frac{(x-y)^2}
    {2{\delta(x)}^2}\Big) \times	\sum_{j=0}^{\theta^2 d(x)N} 
    \frac{(-1)^j}{4^j j!}H_{2j}
    \Bigg(\frac{x-y}{\sqrt{2} \delta(x)}\Bigg),
  \end{equation}
\end{widetext}
%\pagebreak
and the Fourier-space filtering operation \eqref{eq:molli:filtering}
becomes the real-space \textit{mollification} operation given
by
\begin{equation}
  \label{eq:mollification:operation}
    \Phi^d_{p_N}*(S_N[f])(x) = \int \Phi_{p_N}^{d}(x,y) S_N[f](y) dy
\end{equation}
for all points $x$ in the domain \cite{Tanner2006}.\footnote{Within an
  interval $(x-d(x),x+d(x))$, the mollifier \eqref{eq:the:mollifier}
  is symmetric in $x$ and $y$. However, when $x$ and $y$ are
  sufficiently far apart, this symmetry is broken.} The pointwise
convergence of spectral reconstruction $\Phi^d_{p_N}*(S_N[f])(x)$ is
sub-linear near a discontinuity, but improves to exponential
convergence far from it, in accordance with Equation 
\eqref{eq:mollify:convergence} below:
\begin{equation}
    \label{eq:mollify:convergence}
    \lvert \Phi^d_{p_N}*(S_N[f])(x) - S_N[f](x) \rvert \lesssim e ^{- \eta d N}   
\end{equation}
where the constant $\eta$ is dictated by the specific piecewise 
analyticity properties of $f$. Close to the discontinuity $d$ is 
very small and hence the error is large. Far away from the 
discontinuity, however, $d$ increases, forcing the error to become
small.

Note that there is a \textit{different} integral and a
\textit{different} mollification operation for each point $x$. In this
way, mollification is analogous to solving an inverse-problem via
Green's functions methods. We can see this by explicitly rewriting
equation \eqref{eq:mollification:operation} as
\begin{equation}
  \label{eq:mollification:greens:function}
  \Phi_{p,\delta}*(S_N[f])(x) = \int G(x,y) S_N[f](y)dy
\end{equation}
where
% \TODO{Roland suggests ``maybe include `j' in the labels for
%   $\Phi$ so that it is clear that this mollifier is only for a single
%   interval and not all of them.'' I don't understand what he means by
%   this. -JMM}
\begin{equation}
  \label{eq:mollifier:equals:greens}
    G(x,y) = \Phi^d_{p_N}(x,y).
  \end{equation}
  
  For a full formal derivation, we refer the interested Reader to
  Ref. \revone{\cite{Tanner2006, Tadmor2007}} and the references
  therein. For the estimates of convergence of $\Phi^{d}_{p_N}(x,y)$
  see \cite{Tadmor2007} p. 373-374.

\begin{figure}[t]
  \centering
  \includegraphics[width=0.5\textwidth]{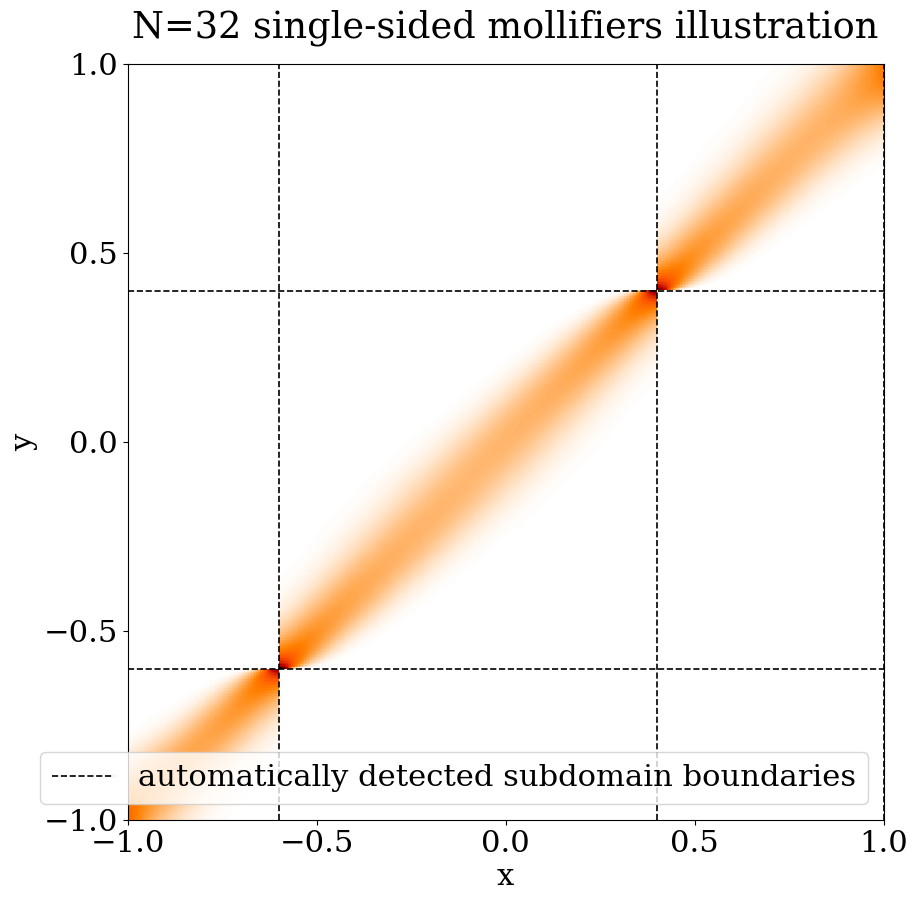}
  \caption{ Illustration of the shape evolution of adaptive single-sided
    mollifiers as they approach a boundary for piecewise
    mollification. Figure presents 500 mollifiers constructed for 
    a~grid of 500 points in the $[-1,1]$ domain. Each point on the
    $x$-axis represents the central point of a mollifier for a
    function with discontinuities at $x=-0.6$ and $x=0.4$. The
    $y$-axis is a copy of the $x$-axis and the color is the amplitude
    of the mollifier. Thus the extent of the color in the $y$-axis is
    the extent of the mollifier in the physical domain. The
    mollifier's support evolves in space, becoming most compact close
    to the function discontinuity. In order to avoid combining
    information from two regions of smoothness across the
    discontinuity, the mollifying kernel is non-zero only on the side
    of discontinuity of interest.}
  \label{fig:mollifiers2D}
\end{figure}

\subsubsection{Preserving discontinuities}
\label{sec:adhoc}

In the limit of $N \longrightarrow \infty$ the dilation operation
defined in Equation \eqref{eq:dilation} becomes a Dirac delta function. However,
for finite $N$, the mollifier is supported on the whole real
line. This means that, even though it becomes narrower near 
discontinuities, it still admixes information across them, 
which leads to and overall smoothing of the edge. 
This is a milder version of the effect
discussed in Section \ref{sec:intro:mollification:filters}.

We would like our method to preserve the discontinuous character of
$f$, at the same time making use of the excellent smoothing properties
of adaptive mollifiers. For this reason we utilize knowledge of jump
locations $\{c_j\}$ and force the mollifier amplitude to vanish
outside the region of smoothness so that:
\begin{equation}
    \label{eq:adhoc}
    \tilde{\Phi}^{d}_{p_N}(x,y):=
    \begin{cases}
        \frac{1}{a(x)}\Phi^{d}_{p_N}(x,y), & 
        \text{if $c_j \leq y \leq c_{j+1}$ }, \\
        0, & \text{otherwise,}
    \end{cases}
\end{equation}
where
\begin{equation}
    a(x) = \int^{c_{j+1}}_{c_j}\Phi^{d}_{p_N}(x,y) dy,
\end{equation}
for a mollifier centered at point $x$ surrounded by edge positions at
$c_j$ and $c_{j+1}$. The parameter $a(x)$ ensures the mollifier is
appropriately normalized to have unit mass.

This approach is \textit{ad hoc}---we do not have a formal
justification for truncating the mollifiers in this way. However, we
believe that the results presented below justify our choice. A formal
analysis of mollifiers which vanish outside the region of smoothness
may reveal a more optimal construction. This is a subject for future
study.

This method results in recovery of truly discontinuous functions,
cured of spurious oscillations.  Illustration of the single-sided
mollifiers is presented in Figures \ref{fig:mollifiers} and
\ref{fig:mollifiers2D}, which show how $\Phi^d_{p_N}(x,y)$ are
forced to zero immediately beyond the theoretical position of the
discontinuity.

\section{Chebyshev Polynomials of First Kind}
\label{sec:decomposition}

In the remainder of this work, we focus on Chebyshev pseudospectral
methods. We therefore spend a few moments to describe them.
%\TODO{We
%  could also move this section to an appendix. That might be
%  appropriate. -JMM}

\subsection{Initial decomposition}
\label{sec:decomposition:initial}

\begin{figure*}
	\centering
	\subfigure[]{
		\label{fig:convergence:gauss}
		\includegraphics[width=.48\textwidth]{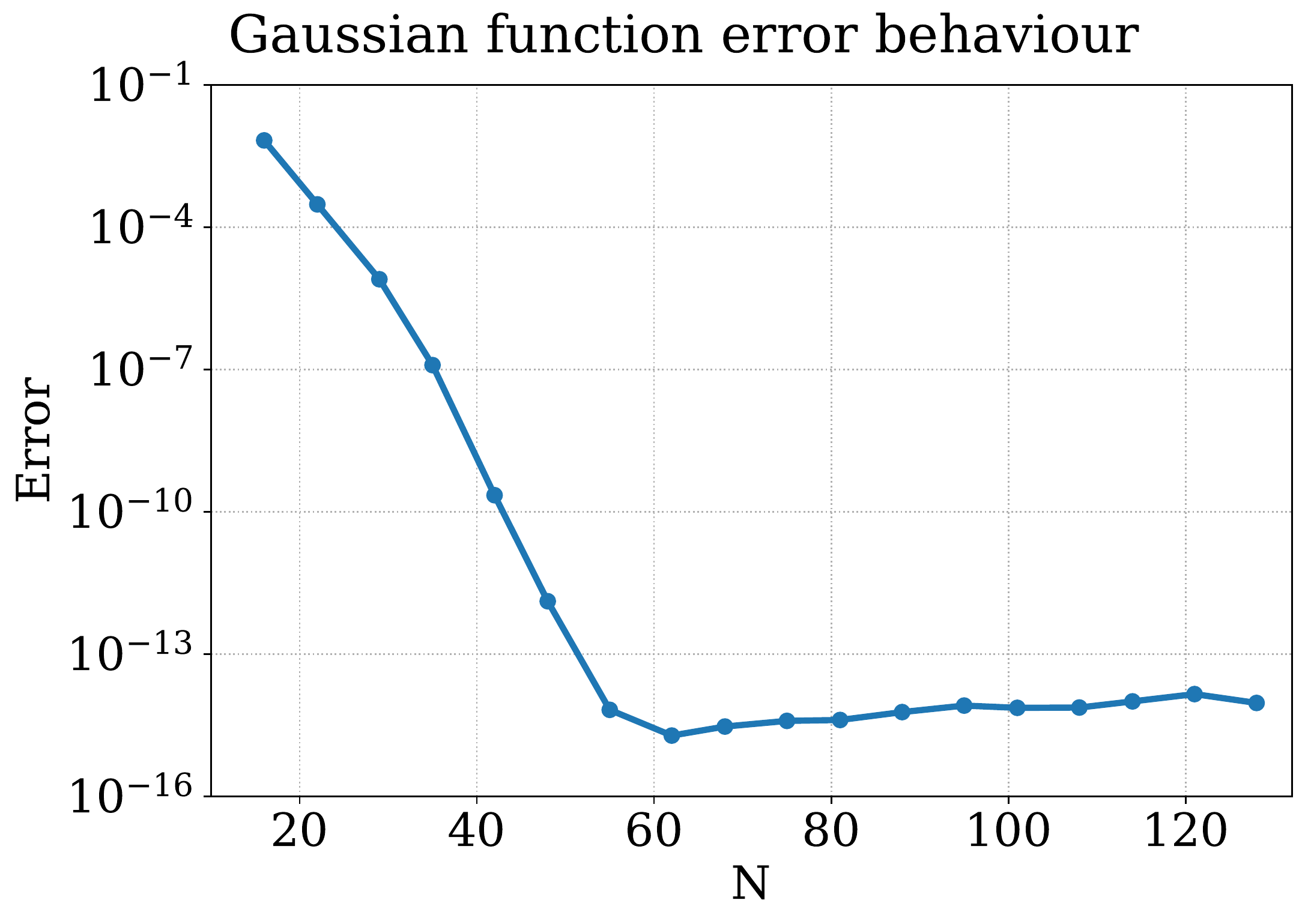}
	}
	\subfigure[]{
		\label{fig:convergence:tophat}   
		\includegraphics[width=.48\textwidth]{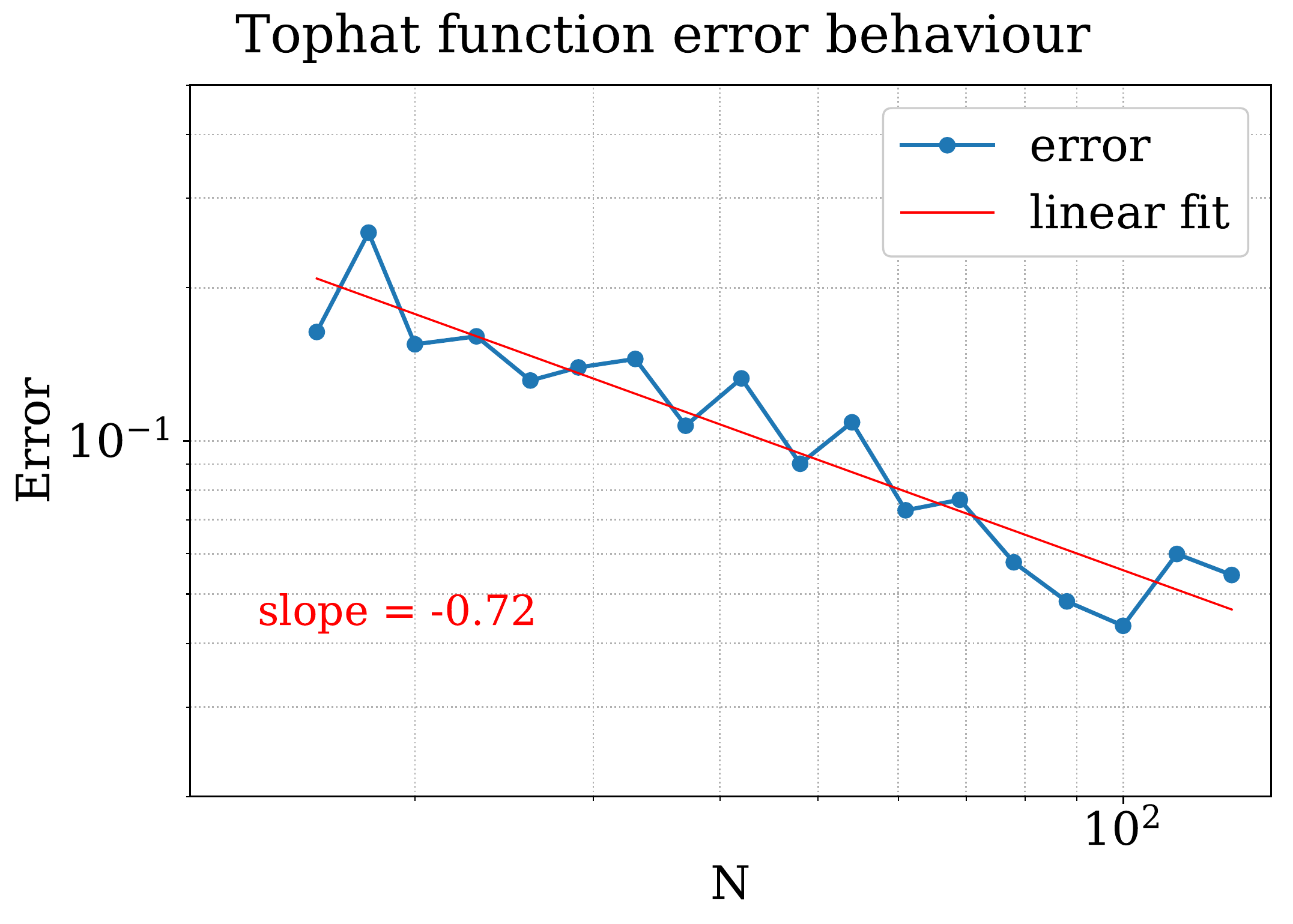}
	}
	\caption{$L_2$-norm difference between the spectral
          reconstruction and true solution as a function of highest
          order of expansion $N$. The norm has been computed using
          high precision quadrature
          integration. \subref{fig:convergence:gauss} shows the smooth
          function case, where the difference decreases exponentially
          with increasing $N$, leveling off at machine precision.
          \subref{fig:convergence:tophat} exhibits slower convergence
          due to presence of discontinuities, which retain
          $\mathcal{O}(1)$ error at the discontinuity for all $N$ and
          allow for $\mathcal{O} (\frac{1}{N})$ convergence away from
          it. The computed slope of linear fit to the plot provides
          the approximate order of convergence.}
	\label{fig:convergence}
\end{figure*}

Chebyshev polynomials of first kind $T_n(x)$ form an orthogonal set
on $[-1,1]$ under the weight $w(x)$ such that:
\begin{equation}
    \label{eq:cheb:ortho}
    \begin{split}
        & \int_{-1}^{1} T_n(x)T_m(x)w(x)dx = 
        \frac{\pi}{2}(1+\delta_{0n})\delta_{nm}, \\ 
        & w(x)=\frac{1}{\sqrt{1-x^2}}.
    \end{split}
\end{equation}

Under restriction to the domain $x\in [-1,1]$ and with the application
of Equation \eqref{eq:cheb:ortho}, Chebyshev polynomials can be used
to approximately represent systems characterized by non-periodic
boundary conditions. (To this purpose the domain must be appropriately 
shifted and rescaled). As discussed in Section
\ref{sec:intro:projection}, such a representation (for smooth
problems) is subject to an error which decays faster than any power of
$N$, as illustrated in Figure \ref{fig:convergence:gauss}. The plot
shows the $L_2$ norm of the difference between the spectral
reconstruction and true solution as a function of the highest order of
expansion $N$ for a Gaussian function $g(x)$ defined below:
\begin{equation}
  \label{eq:def:decay}
    g(x) = e^{\frac{-x^2}{2\sigma^2}},
\end{equation}
where $\sigma=\frac{1}{6}$ to ensure that the relevant piece of $g(x)$
lies well within the domain. The $L_2$ norm decreases
exponentially with increasing $N$, leveling off at $\sim 10^{-15}$,
reaching machine precision. In the case of a discontinuous function,
\begin{equation}
  \label{eq:def:tophat}
  \tau(x) = \begin{cases}
    1, &\text{if } -0.5 \leq x \leq 0.5,\\
    0, &\text{otherwise},
    \end{cases}
\end{equation}
the $L_2$ norm of error behaves differently, decreasing non-uniformly
at a rate $\Ord{N^m}$, where $m\approx-0.7$, as shown in 
Figure~\ref{fig:convergence:tophat}. This trend is consistent with the
recognized influence of discontinuities mentioned in Section
\ref{sec:intro:gibbs}. The Gibbs oscillations present around
discontinuities do not decrease in amplitude, as $N$ is
increased. They persist for any choice of $N$, conserving their
amplitude as they move closer towards the point of jump in the
function.

\subsection{Gauss quadrature and projection}
\label{sec:decomposition:quad}

\begin{figure*}
	\centering
	\includegraphics[width=\textwidth]{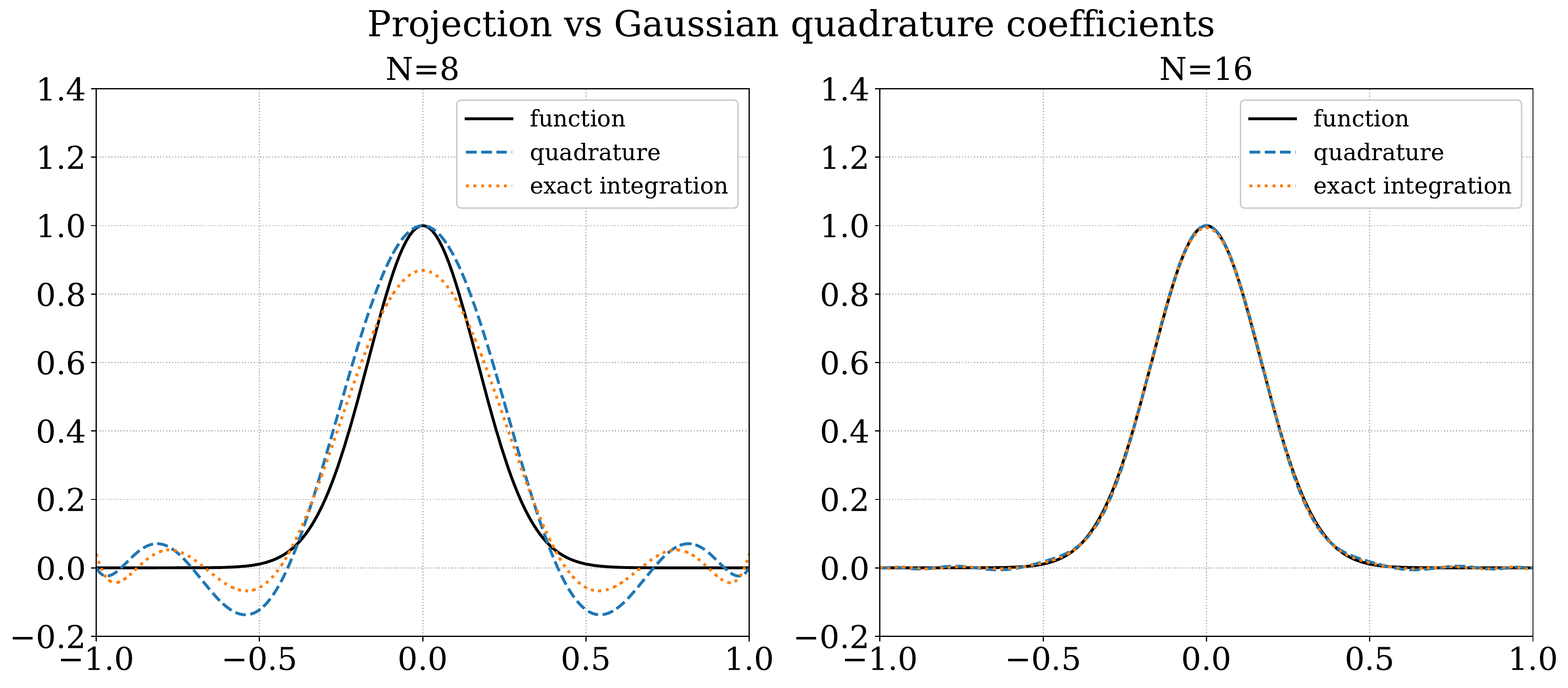}
	\caption{Spectral reconstructions obtained both through
          projection and using Gaussian quadrature. Projection
          coefficients were calculated using high accuracy
          quadrature method in Mathematica
          \cite{Mathematica}, by computing the solution to equation
          \eqref{eq:f_n}. In the smooth function case results differ
          for low orders of $N$, however both reconstructions converge
          onto each other and the true solution exponentially with
          increasing $N$.}
	\label{fig:integral:gauss}
\end{figure*}

\begin{figure*}
	\includegraphics[width=\textwidth]{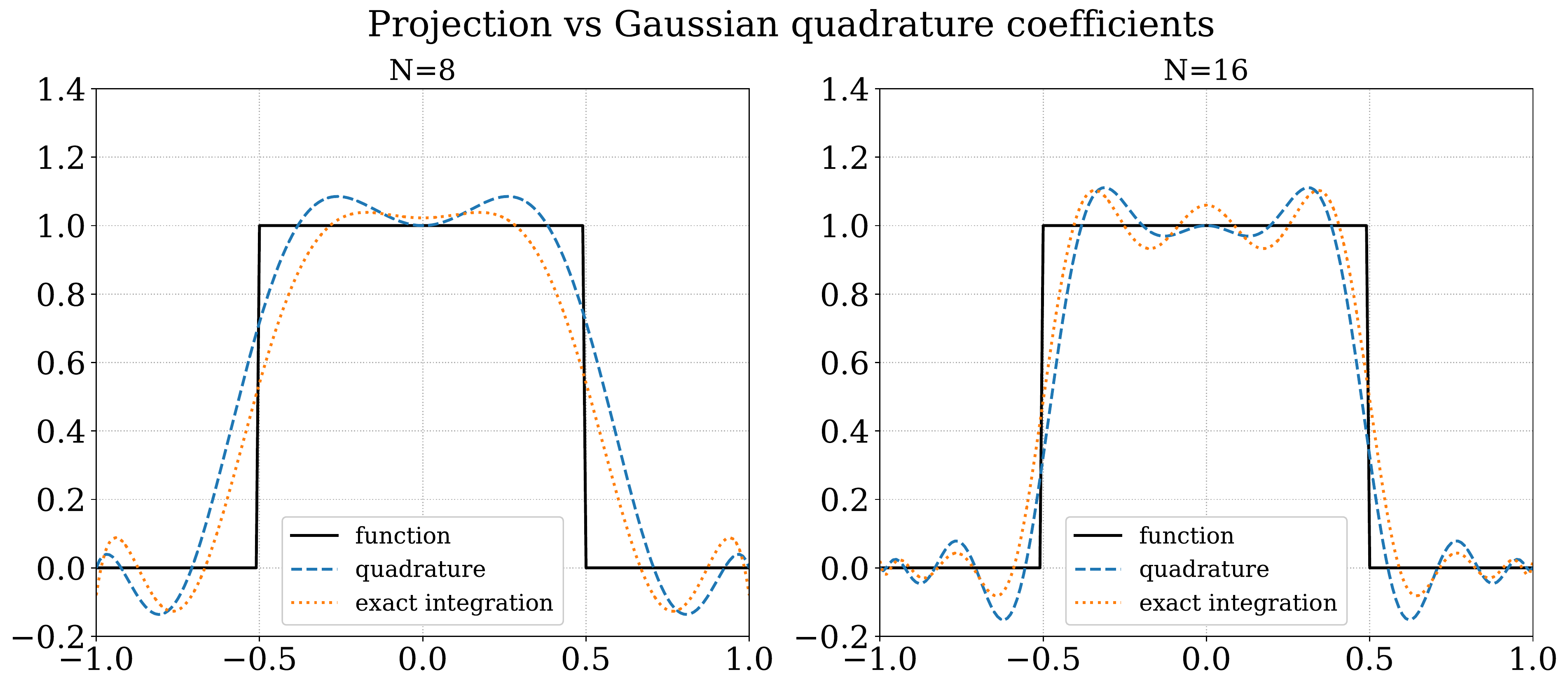}
	\caption{Spectral reconstructions obtained both through
          projection and using Gaussian quadrature. Projection
          coefficients were calculated using high accuracy
          quadrature method in Mathematica
          \cite{Mathematica}, by computing the solution to equation
          \eqref{eq:f_n}.  In the case of a discontinuous function
          reconstructions converge onto each other slower than in
          \ref{fig:integral:gauss}. The high frequency oscillations
          close to the discontinuity never vanish due to the Gibbs
          phenomenon.}
	\label{fig:integral:tophat}
\end{figure*}

The Chebyshev polynomials possess another favorable property, which
facilitates their application in pseudospectral methods introduced in
\ref{sec:intro:projection}. The weights and collocation points
associated with all types of Chebyshev-Gauss quadratures are
completely analytic and can be straightforwardly computed. In this
work, we are using the \textit{Chebyshev-Gauss-Lobatto} quadrature,
which includes collocation points at the domain boundaries. Their
locations and associated weights are given by:
\begin{equation}
  x_i = \cos\frac{\pi i}{N},
  \qquad w_0=w_N=\frac{\pi}{2N}, 
  \qquad w_i=\frac{\pi}{N}.
\end{equation}
As described in Section \ref{sec:intro:projection}, precise
calculation of integrals in order to compute the expansion coefficients is
computationally expensive. The quadrature solution provides convenient
means of handling the operation using Equation \eqref{V}, at the cost
of introducing truncation error in the approximation of the
integral. In the case of smooth functions the pseudospectral and
spectral reconstructions converge to the true solution at an
exponential rate. In other words, they converge to each other, as
illustrated in Figure \ref{fig:integral:gauss}.

Figure \ref{fig:integral:gauss} presents spectral reconstructions of a
smooth Gaussian function using two sets of expansion coefficients. The
quadrature ones are computed using the RHS of Equation
\eqref{gauss_quad}, while the exact integration ones are computed in
Mathematica \cite{Mathematica} as an inner product between the
Gaussian and a given Chebyshev polynomial, using its high accuracy
quadrature method. The same procedure is repeated in the case of a
discontinuous function in Figure \ref{fig:integral:tophat}. Because of
the Gibbs phenomenon, neither the spectral nor the pseudospectral
projection converges onto the original function across the entire
domain and the projections differ from each other. This effect needs
to be recognized as an additional source of error in any method
applied to discontinuous functions.

\subsection{Spatial derivatives}
\label{sec:decomposition:deriv}

\begin{figure*}
	\subfigure[]{
		\label{fig:deriv:gauss}
		\includegraphics[width=0.48\textwidth]{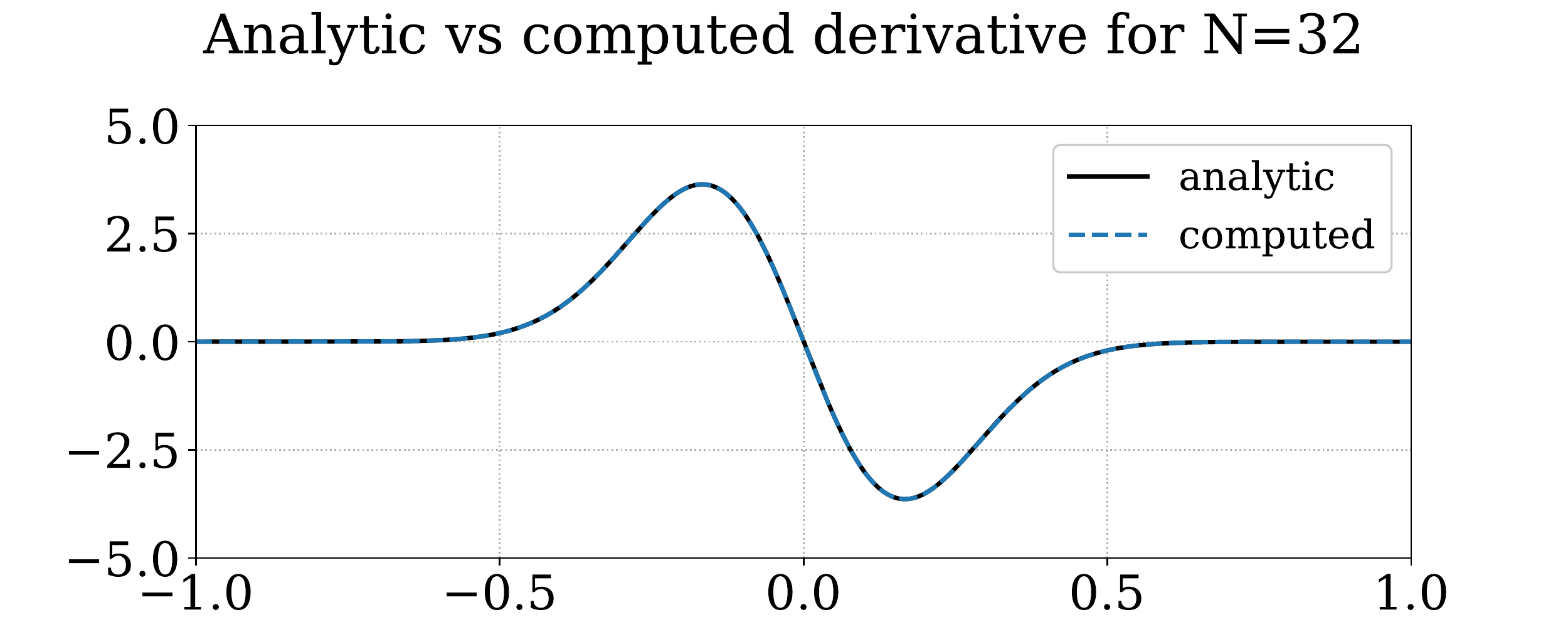}
	}
	\subfigure[]{
		\label{fig:deriv:tophat}
		\includegraphics[width=0.48\textwidth]{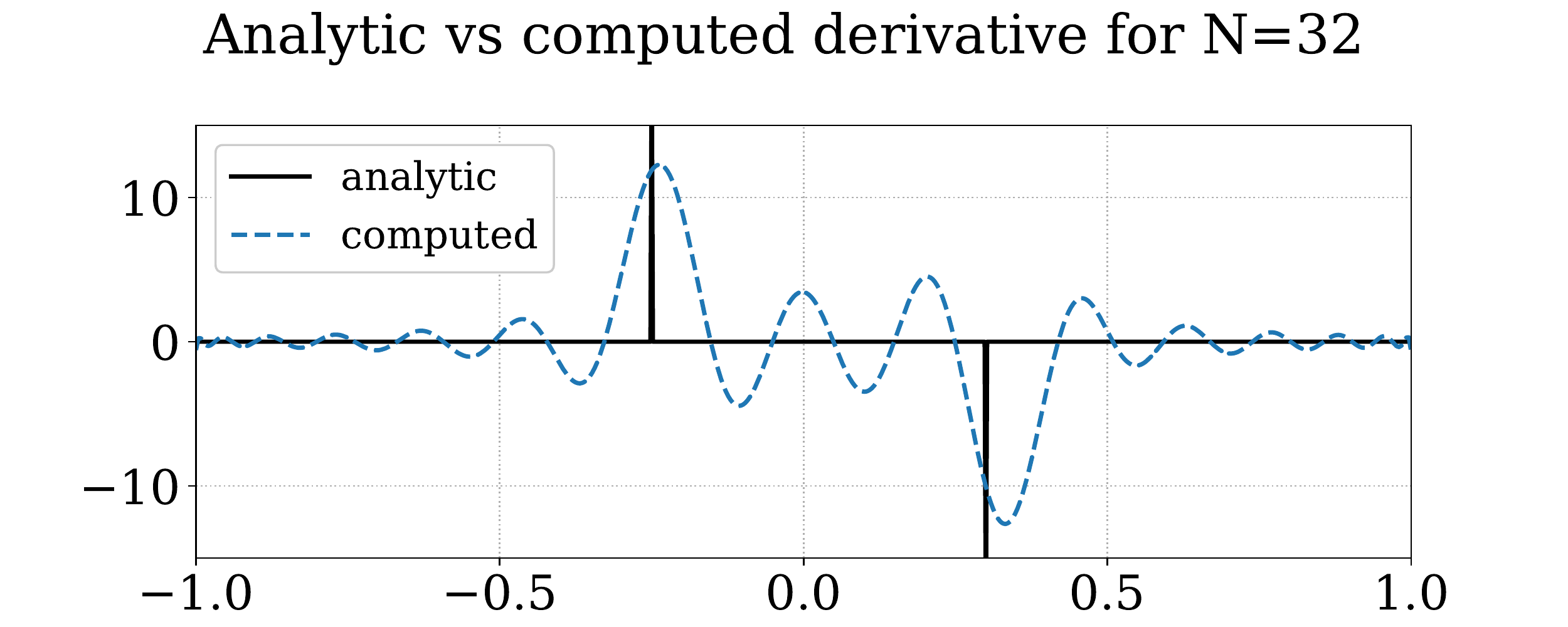}
	}
	\caption{
          Visual comparison between analytic derivative of function defined
          in Equation \eqref{eq:def:decay} and one
          computed using Equation \eqref{eq:D:modal} for a smooth and
          discontinuous function case. \subref{fig:deriv:gauss} shows
          perfect overlap between the analytic and computed result for
          $N=32$ highest order of expansion. \subref{fig:deriv:tophat}
          illustrates the influence of the Gibbs phenomenon on the
          derivative of a~discontinuous function. Instead of obtaining
          two delta functions at the locations of the discontinuities,
          one observes an oscillating pattern with peaks of highest
          amplitude located close to the discontinuity. Note that the
          peaks are not located exactly at the discontinuity. There is
          some ``phase error.''}
	\label{fig:deriv}
\end{figure*}

Spectral reconstructions of derivatives of smooth functions
inherit the supreme error convergence properties, visualized in Figure
\ref{fig:deriv:gauss}. The result computed using Equation
\eqref{eq:D:modal} completely overlaps with the analytic result,
with the $L_2$-norm of the difference of order $\Ord{10^{-5}}$,
as expected given Figure \ref{fig:convergence:gauss}.

With discontinuous functions, however, the Gibbs phenomenon
influences the derivative across the entire domain, producing
spurious oscillations. This effect is presented in Figure
\ref{fig:deriv:tophat}, where the computed spectral reconstruction
differs significantly from the analytic solution.  It is also
important to note that the information about the exact position of the
discontinuity is subject to uncertainty due to the significant size of
spacing between collocation points. This effect is less pronounced
closer to the domain boundaries, where collocation points are more
closely spaced.

\section{Edge Detection}
\label{sec:edge:detection}

As we saw in Section \ref{sec:decomposition:deriv}, the spectral
derivative of a~discontinuous function is not well behaved. This is
unfortunate, since the (infinitely) large gradient at a discontinuity
provides an excellent means of detecting one. Here we describe the
technique, first developed by Gelb
\cite{Gelb1999a,Gelb2005a,gelb2006robust}, for regularizing these
spectral derivatives and using them as a way to localize
discontinuities in spectral data.

\subsection{Method}
\label{sec:edge:detection:method}

We define the jump function $[f](x)$ of a piecewise smooth function 
$f(\cdot)$ to be:
\begin{equation}
    [f](x):= f(x^{+}) - f(x^{-}).
\end{equation}
In \cite{Gelb2001a}, Gelb showed that for admissible
\textit{concentration factors} $\mu(\frac{k}{N})$ and spectral
coefficients $\hat{f}_k$, in the basis of polynomial derivatives
$T^{\prime}_k(x)$, the expression below converges pointwise to
$[f](x)$:
\begin{equation}
    \label{eq:edge:concentration}
    \frac{\pi\sqrt{1-x^2}}{N}\sum^N_{k=1} \mu\Bigg(\frac{k}{N}\Bigg)
    \hat{f}(k){T^\prime}_k(x) \longrightarrow [f](x).
\end{equation}
The basic idea here is to filter (\`a la Vandeven \cite{Vandeven1991a})
the spectral expansion of the gradient of $f$ so that the spurious
oscillations in the gradient disappear.

As discussed in \cite{Gelb1999a, Gelb2001a, Gelb2008a} there exist
multiple concentration factors with convergence rate of $\sim
\mathcal{O}(\frac{\log{N}}{N})$, which exhibit 
different behaviour close to the discontinuity. Each
$\mu(\frac{k}{N})$ produces its own oscillatory pattern near the jump
discontinuities, often trading exponential convergence away from the
jump for decreased oscillations close to it. In order to exploit
both good convergence and no-oscillation features of individual
factors, Gelb and collaborators \cite{Gelb2008a} designed a
\textit{minmod function} as
\begin{widetext}
\begin{equation}
    \label{eq:minmod}
    \text{minmod}(f_1,\ldots,f_j)(x):=
    \begin{cases}
        \min(f_1(x),\ldots,f_j(x)), & \text{if $f_1(x),\ldots,f_j(x)>0$,} \\
        \max(f_1(x),\ldots,f_j(x)), & \text{if $f_1(x),\ldots,f_j(x)<0$,} \\
        0, & \text{otherwise,}
    \end{cases}
\end{equation}
\end{widetext}
where each $f_i$ corresponds to the jump function approximation
computed through Equation \eqref{eq:edge:concentration} with a
different concentration factor $\mu$. Note that the minmod function is
a pointwise object, since the chosen $f$ can be different at each point
$x$. The estimated convergence rate of $\text{minmod}(f\ldots)$ given
in Equation \eqref{eq:minmod} to $[f](x)$ is given by the lowest
formal bound on the individual concentration methods and evaluates to
$\sim\mathcal{O}(\frac{\log N}{N})$ (for a formal proof, see
\cite{Gelb2008a} and references therein).

In our investigation, we follow \cite{Gelb2008a} and make use of Equation \eqref{eq:minmod} by
applying concentration factors $\mu(\frac{k}{N})=\mu(\eta)$ first
defined in \cite{Gelb1999a} as
\begin{equation}
    \label{eq:cf:trig}
    \mu_{\text{trig}}(\eta) = \frac{\pi\ \sin (\beta \eta)}{\sin(\beta)},
\end{equation}
\begin{equation}
    \label{eq:cf:poly}
    \mu_{\text{poly}}(\eta) = \pi\ \eta^p,
\end{equation}
and
\begin{equation}
    \label{eq:cf:expon}
    \mu_{\text{exp}}(\eta) = \gamma\eta\ \exp\Bigg(\frac{1}
    {\alpha\eta(\eta-1)}\Bigg),
\end{equation}
where
\begin{equation*}
    \gamma = \frac{\pi}{\int_{\epsilon}^{1-\epsilon} 
    \exp(\frac{1}{\alpha \tau (\tau-1)})d \tau}.
\end{equation*}
Certainly other concentration factors are valid and
  the application of additional factors is worth investigating. As in \cite{Gelb2008a}, we set
$\beta=\pi$, $p=1$, $\alpha=6$ and $\epsilon=\frac{1}{N}$ in our
tests. Results of our computations are discussed in Section
\ref{sec:edge:detection:performance} below.

\begin{figure}
  \centering
  \subfigure[]{
  \label{fig:jump:trig}
  \includegraphics[width=.5\textwidth]{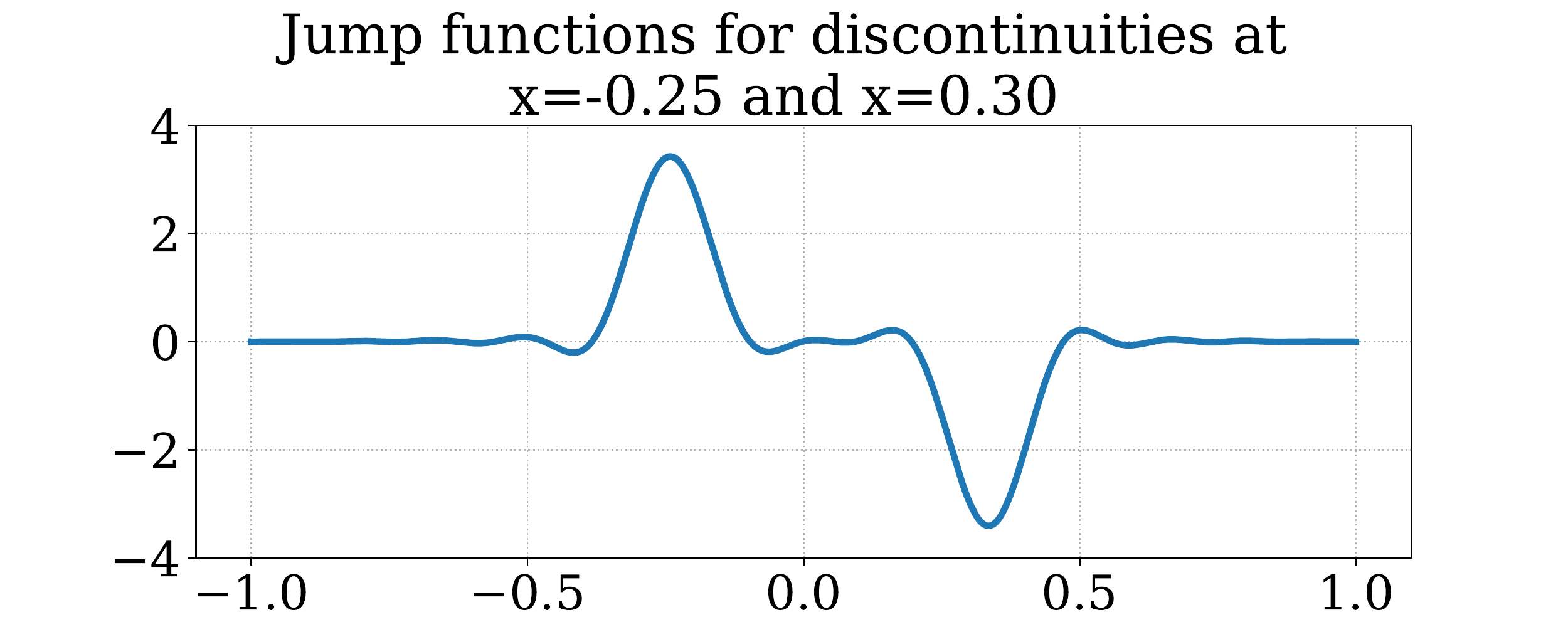}
  }
  \subfigure[]{
  \label{fig:jump:poly}   
   \includegraphics[width=.5\textwidth]{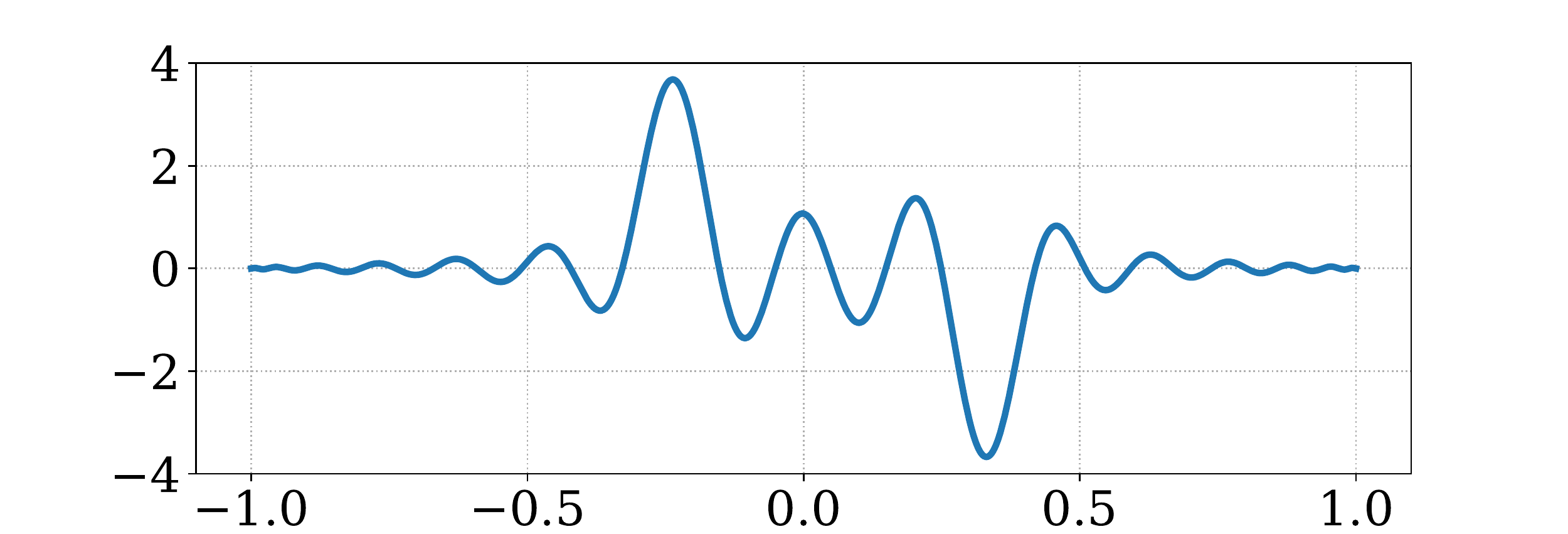}
   }
  \subfigure[]{
  \label{fig:jump:expon}   
   \includegraphics[width=.5\textwidth]{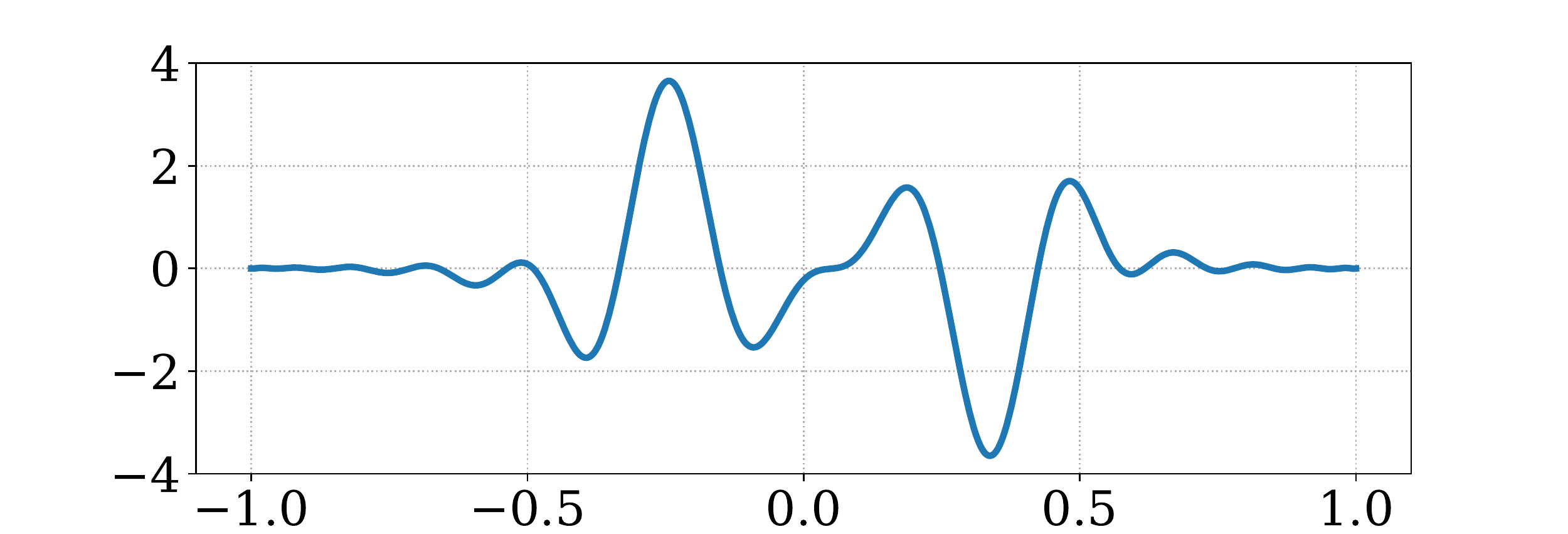}
   }
   \caption{Jump functions computed using Equations
     \eqref{eq:cf:trig}, \eqref{eq:cf:poly}, and \eqref{eq:cf:expon}.
     \subref{fig:jump:trig}\,presents the result of applying
     $\mu_{\text{trig}}$, while \subref{fig:jump:poly} and
     \subref{fig:jump:expon} correspond to $\mu_{\text{poly}}$ and
     $\mu_{\text{exp}}$ respectively. Each jump function is characterised by
     different oscillation patterns, which allow for the powerful
     performance of the minmod recipe in \eqref{eq:minmod}.}
  \label{fig:jump}
\end{figure}

\subsection{Efficacy}
\label{sec:edge:detection:performance}

\begin{figure}
  \centering
  \subfigure[]{
  \label{fig:minmod:good}
  \includegraphics[width=.5\textwidth]{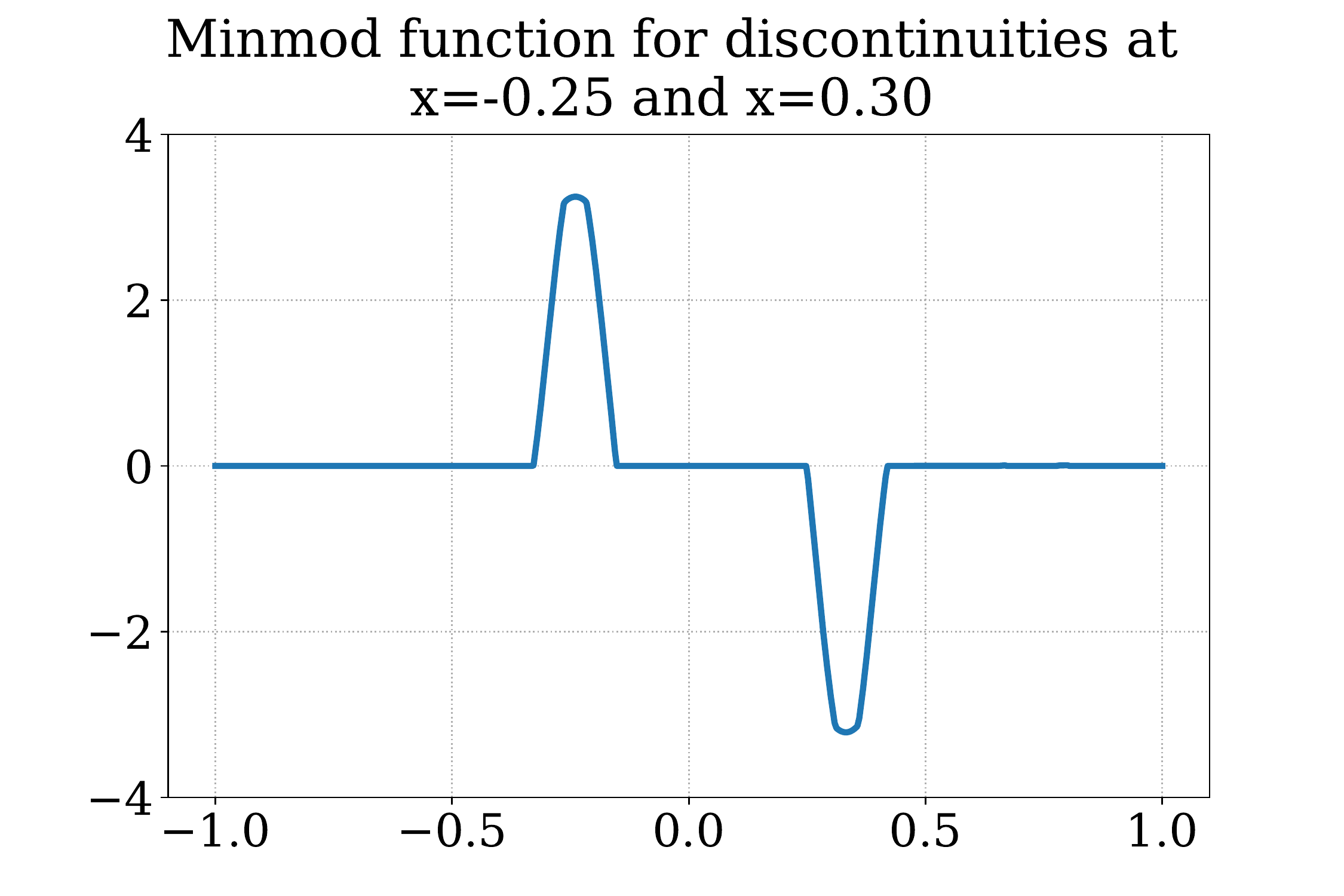}
  }
  \subfigure[]{
  \label{fig:minmod:bad}   
   \includegraphics[width=.5\textwidth]{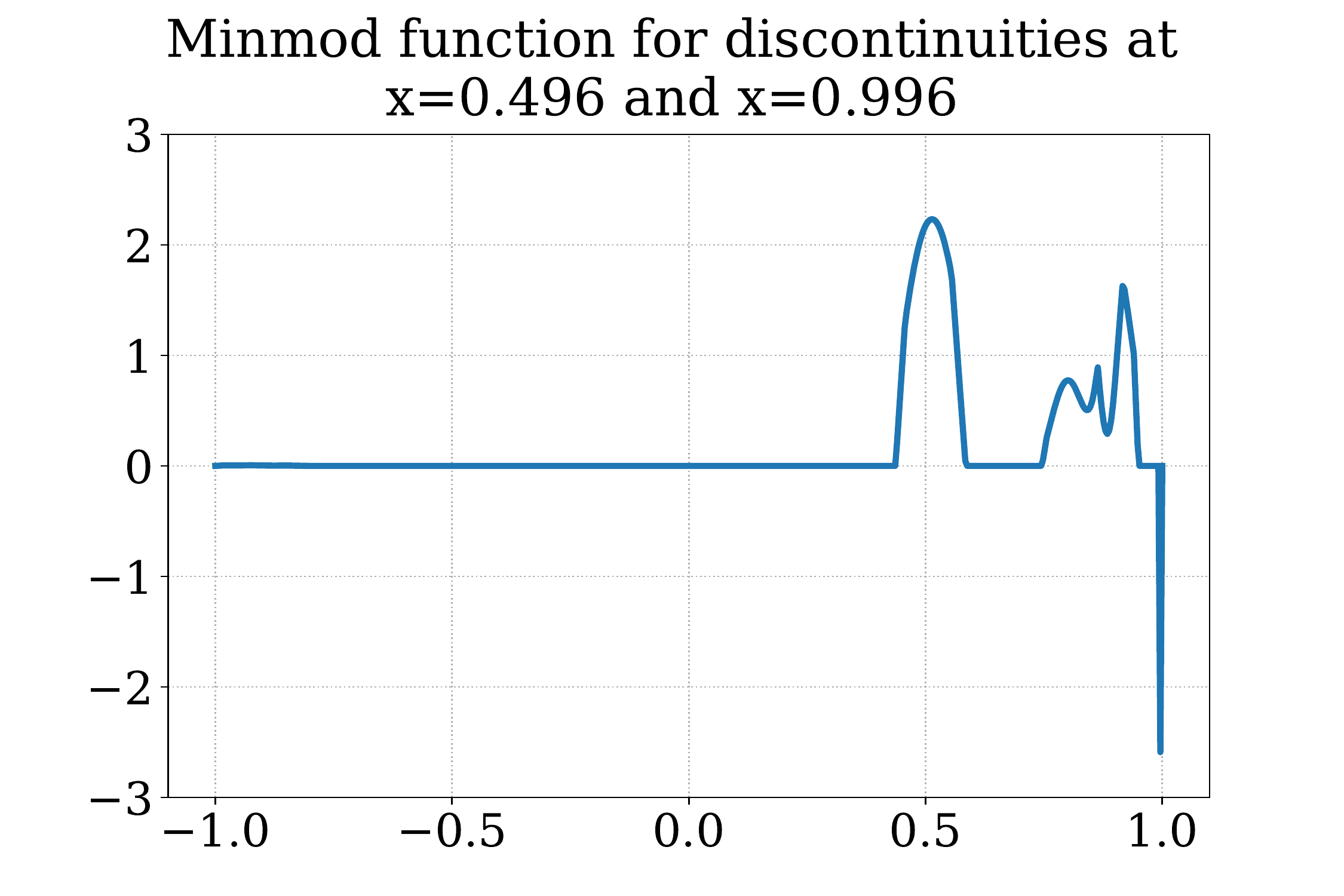}
   }
   \caption{Minmod function computed following \eqref{eq:minmod}, using
     concentration factors in \eqref{eq:cf:trig}, \eqref{eq:cf:poly},
     \eqref{eq:cf:expon}. \subref{fig:minmod:good} presents the case of
     discontinuities located well within the domain, where
     identification of \textit{minmod} function extrema allows for
     straightforward computation of the discontinuity location. Errors
     are quantified in Figure \ref{fig:edgeloc:centre}.
     \subref{fig:minmod:bad} illustrates the inability to easily
     identify the true discontinuity location because of the minmod
     function being polluted by spurious extrema. Errors associated
     with this more difficult case are quantified in Figure
     \ref{fig:edgeloc:boundary}.}
  \label{fig:minmod}
\end{figure}

% \TODO{Add transition paragraph/sentence -JMM}

The concentration factors given in Equations \eqref{eq:cf:trig},
\eqref{eq:cf:poly}, and \eqref{eq:cf:expon} produce jump function
approximations akin to those presented in Figure \ref{fig:jump}. The
output is then fed into Equation \eqref{eq:minmod}, which results in
the functions visualised in Figure \ref{fig:minmod}. Location of
discontinuities is then found by
searching for local extrema in the minmod result.

\begin{figure}
\centering
\includegraphics[width=.5\textwidth]{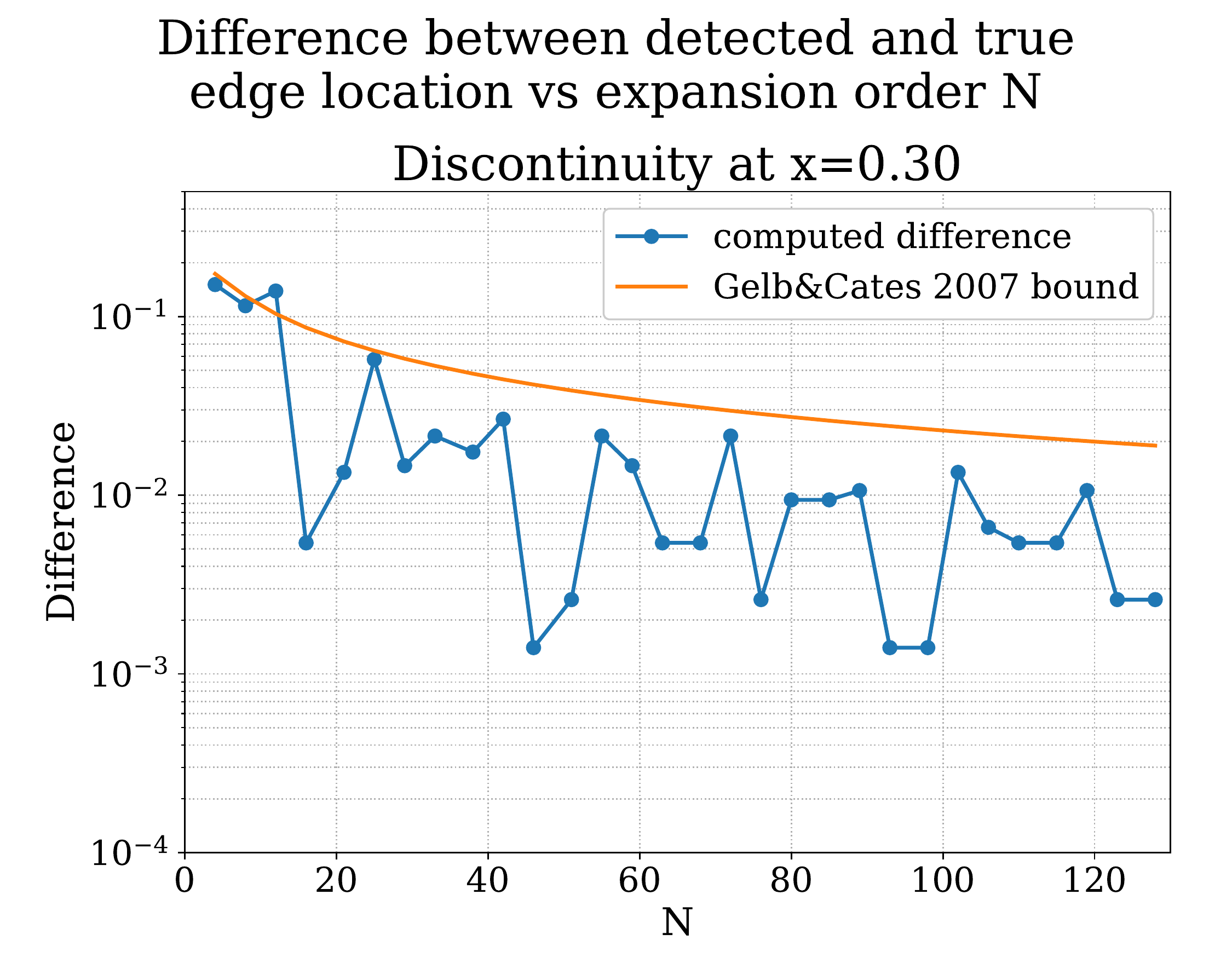}
\caption{Difference between the true discontinuity location and the
  one computed using edge detection recipe provided by
  \cite{Gelb2008a} for a discontinuity located \textit{well within the
    domain}, as a function of order of expansion $N$. There were two
  discontinuities of the same jump magnitude, present at $x=-0.25$ and
  $x=0.30$ on the domain $[-1,1]$. The minmod function is sampled on a
  real space grid of 500 points. The edge detection algorithm performs
  well, with the error convergence upper bound following
  $\mathcal{O}(\frac{\log N}{N})$ relation derived in
  \cite{Gelb2008a}.}
  \label{fig:edgeloc:centre}
\end{figure}

In simple cases, when the discontinuities are located well within the
domain and have comparable magnitudes, the method proves successful,
returning discontinuity position with an error of
$\sim\mathcal{O}(10^{-2})$ for $N=32$ highest order of expansion.  To
illustrate this claim, we present the error in discontinuity location
and its behaviour with increasing $N$ in Figure
\ref{fig:edgeloc:centre}, where the jumps occur at $x=-0.25$ and
$x=0.30$. Oscillations in the error behaviour are due to changing
spatial distribution of collocation points with increasing number of
modes. The envelope, however, seems to obey the
$\sim\mathcal{O}(\frac{\log N}{N})$ error convergence predicted by
theoretical arguments in \cite{Gelb2008a}. As shown in Figure
\ref{fig:minmod:good}, the minmod function cleanly predicts two peaks,
identifying the edges without issue.

\begin{figure}
 \includegraphics[width=.5\textwidth]{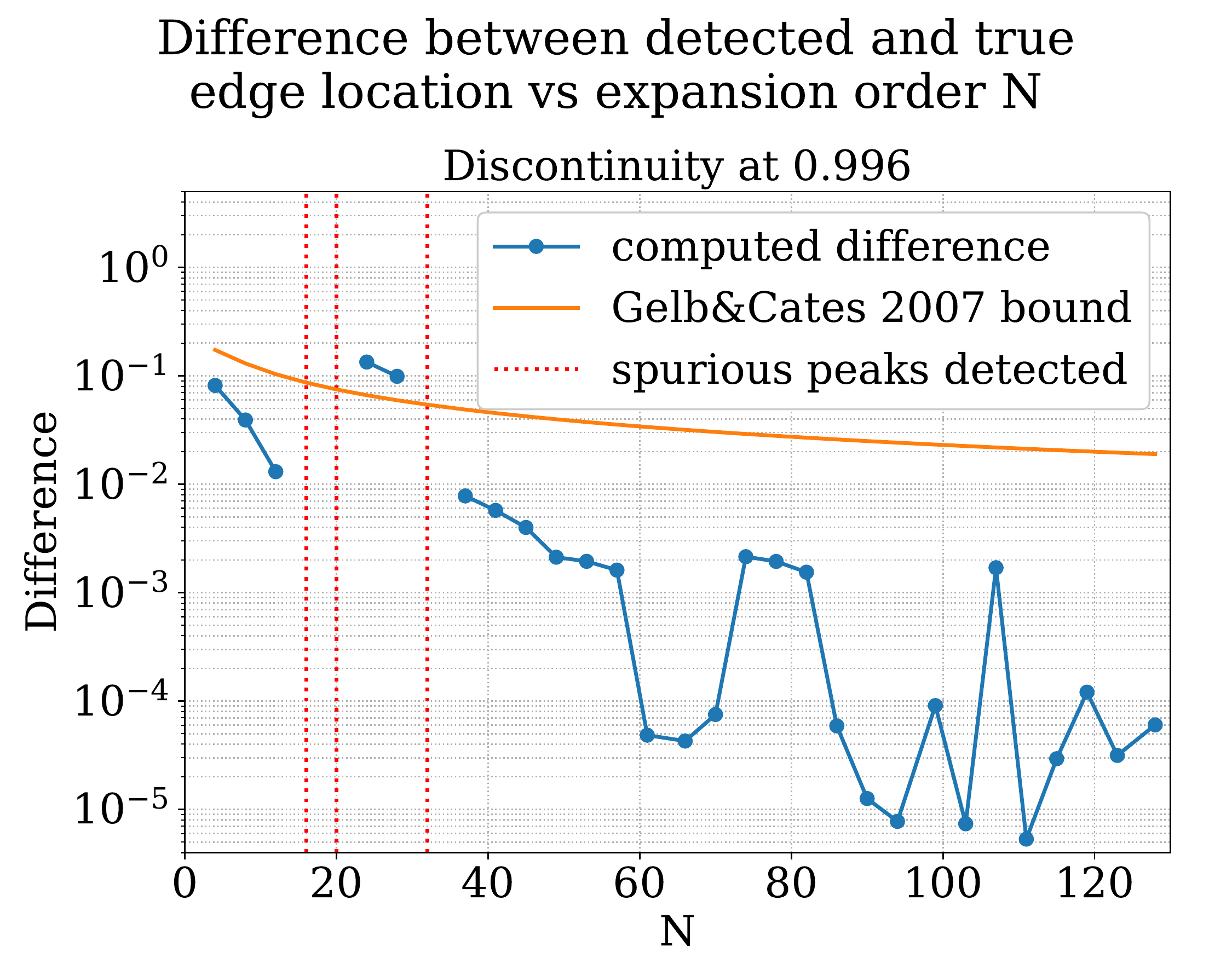}
 \caption{Difference between the true discontinuity location and 
    the one computed using edge detection recipe provided by 
    \cite{Gelb2008a} for a discontinuity located 
    \textit{close to the domain boundary}, as a function of order 
    of expansion $N$. There were two discontinuities of the same
   jump magnitude, present at $x=0.496$ and $x=0.996$ on the domain
   $[-1,1]$. The minmod function was sampled on a real space grid of 500
   points. The edge detection algorithm performs well for high enough
   expansion order $N$ with the absolute error lower than in the case
   of discontinuity location well within the domain. However, it finds
   spurious peaks and fails to find the physical discontinuity for
   $N=16,20,32$.}
  \label{fig:edgeloc:boundary}
\end{figure}

Discontinuity location accuracy is not only a function of the order of
expansion $N$, but also the spatial position of the jump within the
domain. In the center of the domain, where the collocation point
density is smallest, the accuracy is naturally lowest. As the
discontinuity moves closer to the domain boundary, the magnitude of
difference between detected and true jump location decreases due to
higher density of collocation points---compare the absolute values of
the error in Figure \ref{fig:edgeloc:boundary} with those in Figure
\ref{fig:edgeloc:centre}. However, multiple peaks begin to appear in
the minmod function, due to poor behaviour of jump function
approximations close to the domain boundary. These spurious peaks can
make edge detection unreliable near domain boundaries. This issue is
well illustrated in Figure \ref{fig:minmod:bad}, where clear
identification of the correct maximum is not possible. Knowledge of
spectral expansion coefficients at this point is not enough to
identify a single discontinuity in this regime.

In the attempt to mitigate this effect and allow for correct
discontinuity location up to the domain boundaries, we solved the
problem heuristically. Our approach utilized spectral reconstruction
of the derivative of the function of interest, which contains 
information about the expected sign of the jump. Despite the presence
of Gibbs phenomenon, we were able to discard 
spurious minmod peaks by comparing their sign against neighboring 
extrema in the derivative. In the case of fluid shocks,
one might be able to use shock indicators for further guidance. This
approach is not yet robust and requires more experimentation in more
realistic settings. Another way to mitigate this difficulty might be
via nonlinear enhancement, as described in \cite{Gelb2001a}.

These spurious peaks provide further difficulties in the presence of
multiple discontinuities of varying magnitudes. It is not always clear
how to distinguish a smaller, physical peak in the minmod function
from a spurious one. Moreover, a small discontinuity may be washed out
entirely if it is accompanied by a larger one.

\section{Mollification}
\label{sec:mollifiers}

\begin{figure*}
  \centering
  \subfigure[]{
  \label{fig:mollified32:centre}
  \includegraphics[width=.48\textwidth]{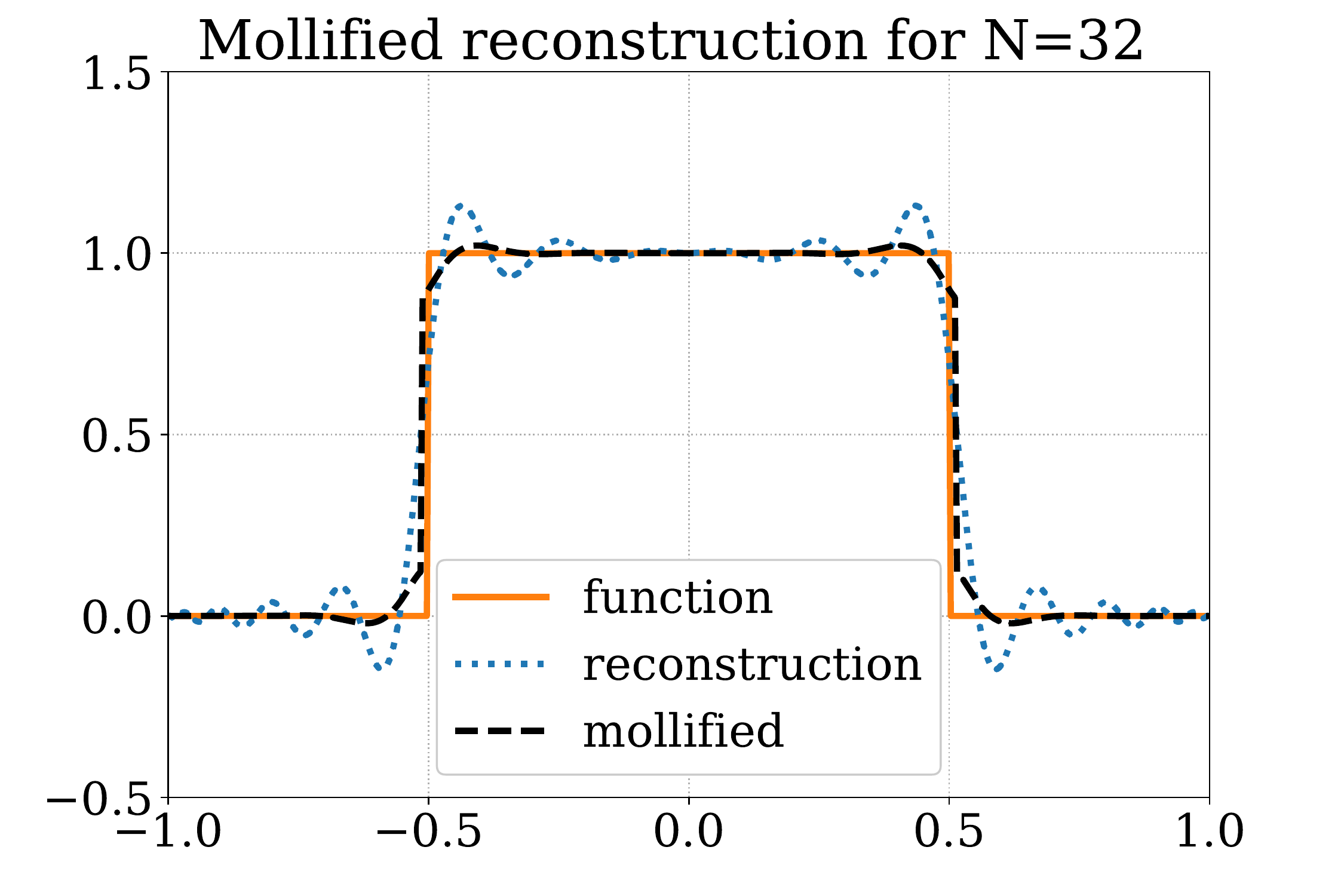}
  }
  \subfigure[]{
  \label{fig:mollified32:edge}   
   \includegraphics[width=.48\textwidth]{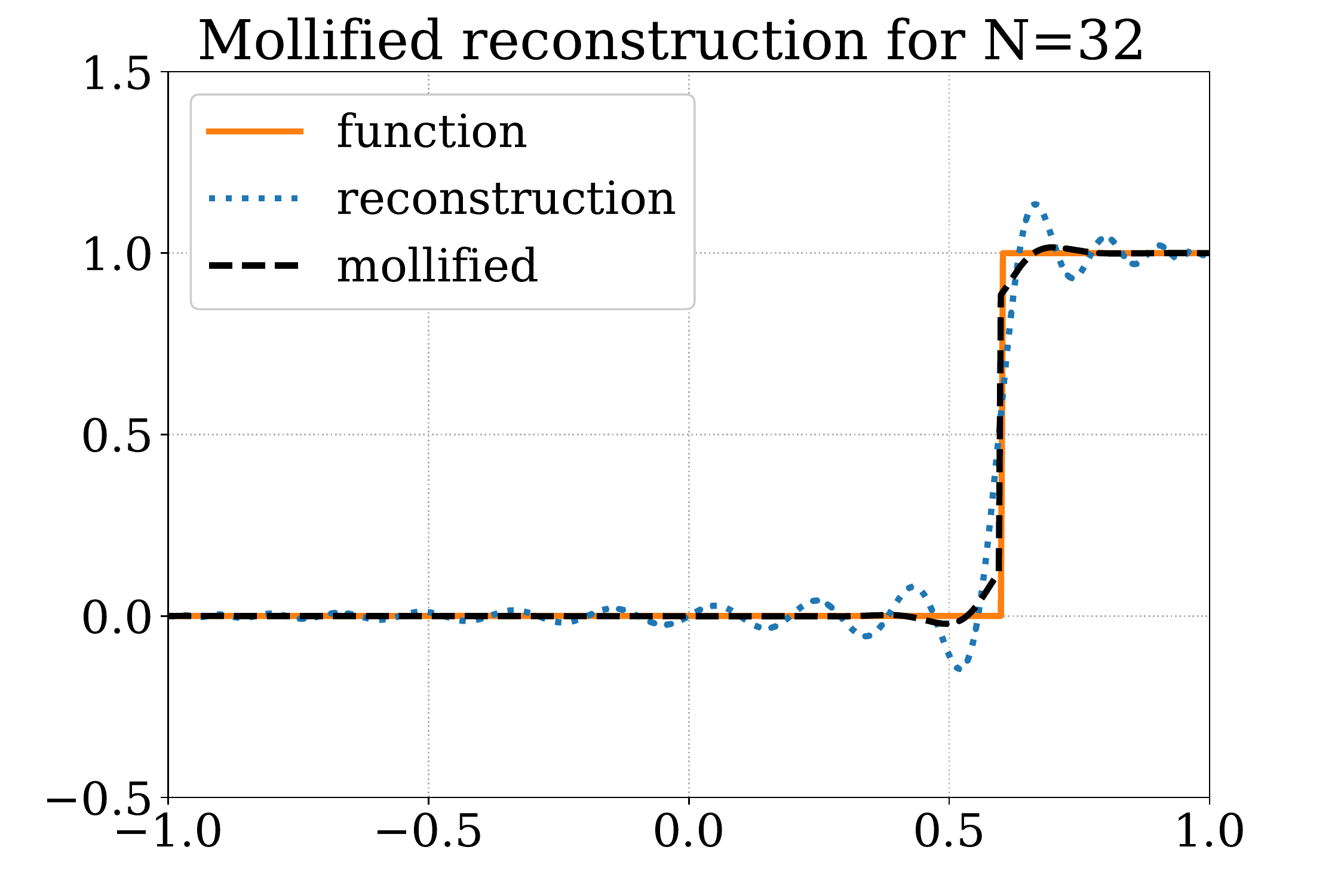}
   }
   \caption{Effect of piecewise mollification applied to a top-hat
     function projected onto a basis of $N=32$ Chebyshev
     polynomials. Piecewise mollification (thick solid) smooths out
     the Gibbs oscillations in the reconstruction (dashed) and
     preserves the discontinuous properties of the original function
     (fine solid), irrespective of a shallow slope of the
     reconstruction. \subref{fig:mollified32:centre} Two
     discontinuities located well within the domain.
     \subref{fig:mollified32:edge} A single discontinuity well within
     the domain.}
  \label{fig:mollified32}
\end{figure*}

\begin{figure*}
  \centering
  \includegraphics[width=\textwidth]{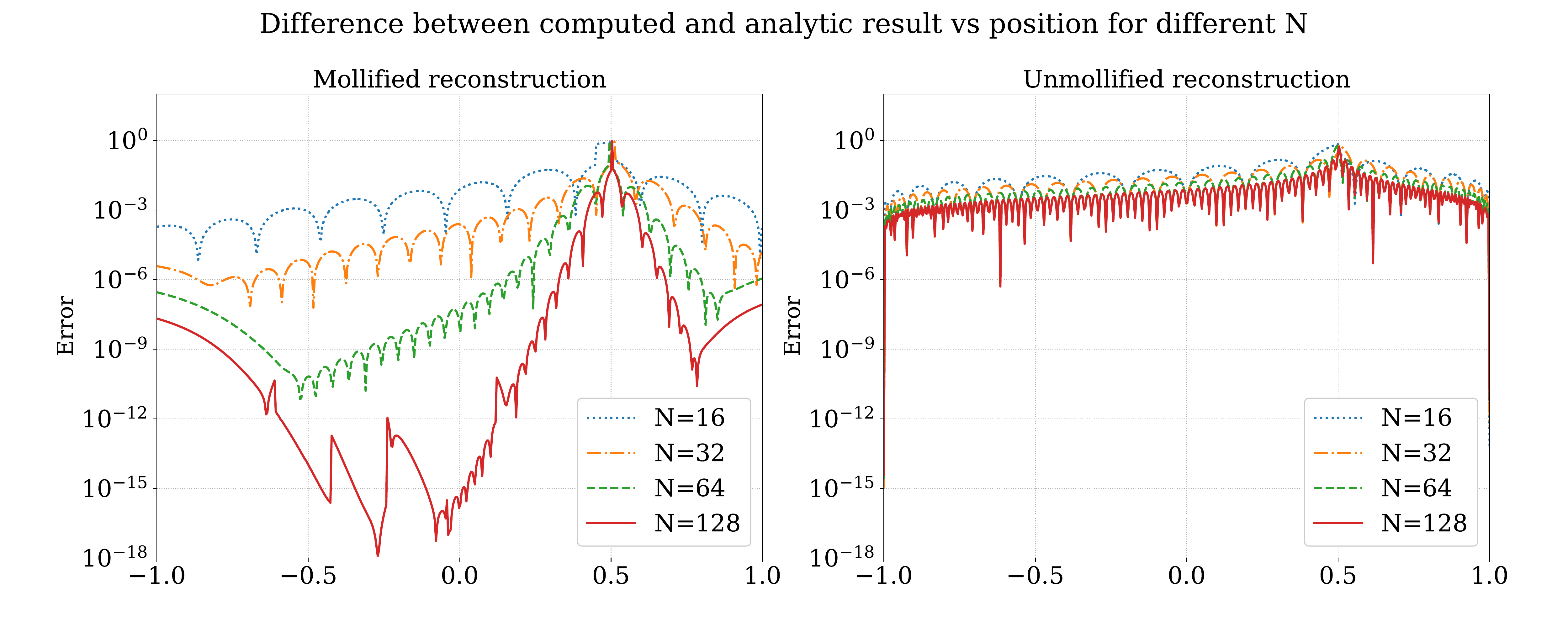}
  \caption{Pointwise difference between the analytic and computed
    result for mollified (left) and unmollified (right) spectral
    reconstruction sampled on an equally spaced grid of 500 points in
    the domain $[-1.1]$. The panels show the case of a single
    discontinuity within the domain, as presented in Figure
    \ref{fig:mollified32:edge}. The error in mollified reconstruction
    converges exponentially away from the discontinuity, as shown in
    Figure \ref{fig:mollifier:convergence}, while the unmollified
    reconstruction suffers poor convergence. Performance of
    mollification with $N=16$ modes involved is comparable to that of
    the $N=128$ unmollified spectral reconstruction.}
    \label{fig:mollified_unmollified}
\end{figure*}

\begin{figure}
	\centering
	\includegraphics[width=0.5\textwidth]{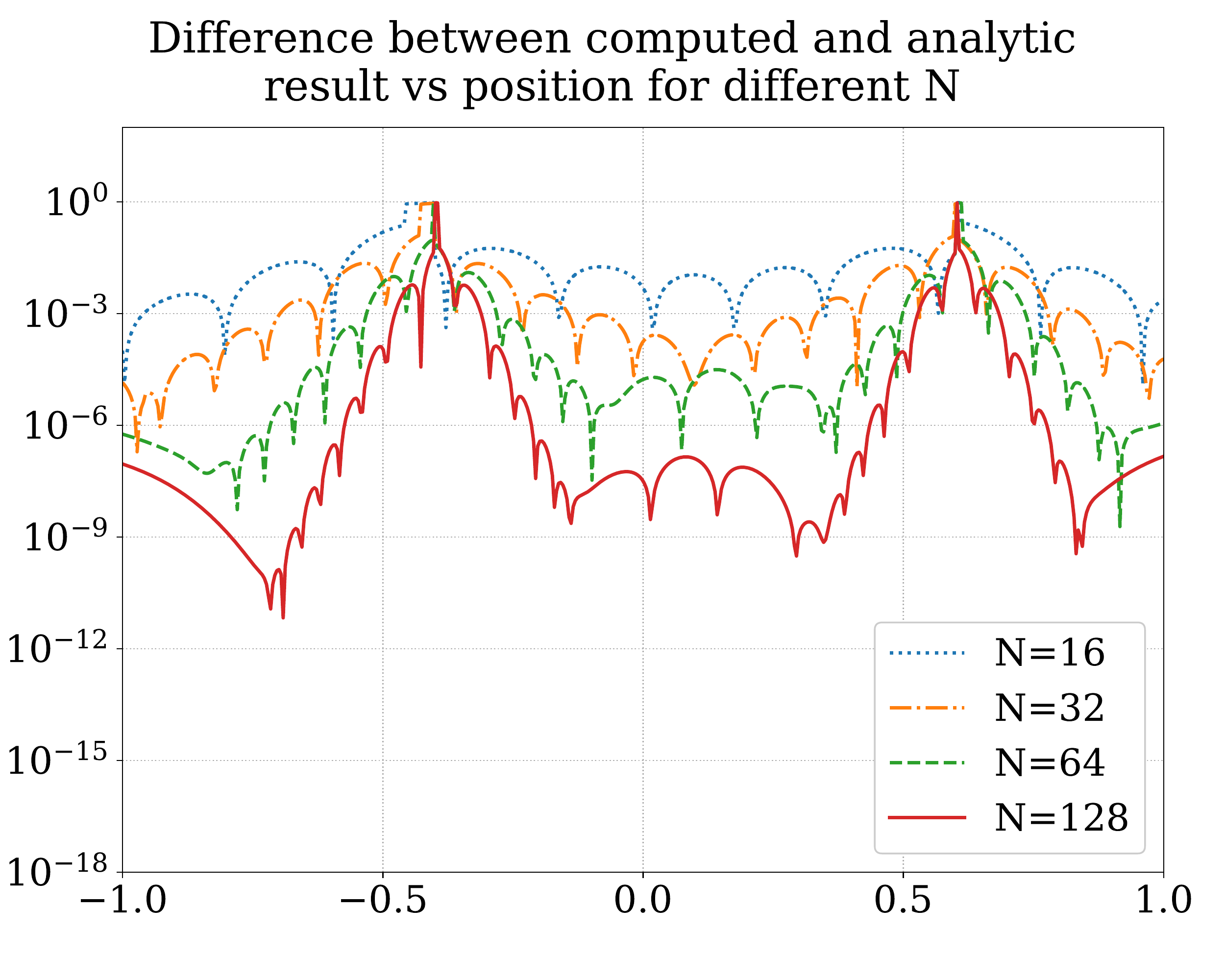}
	\caption{Pointwise difference between the analytic result and
          mollified spectral reconstruction sampled on an equally
          spaced grid of 500 points in the domain $[-1,1]$.  Here we
          have two discontinuities within the domain, as seen in
          Figure \ref{fig:mollified32:centre}. The error exhibits
          $\Ord{N^{-4}}$ decay between the discontinuities and
          $\Ord{N^{-3}}$ decay at domain boundaries, which is further
          demonstrated in Figures \ref{fig:mollifier:convergence} and
          \ref{fig:mollifier:convergence:boundary}.}
	\label{fig:mollified:pointwise:error}
\end{figure}

\begin{figure*}
  \centering
  \includegraphics[width=\textwidth]{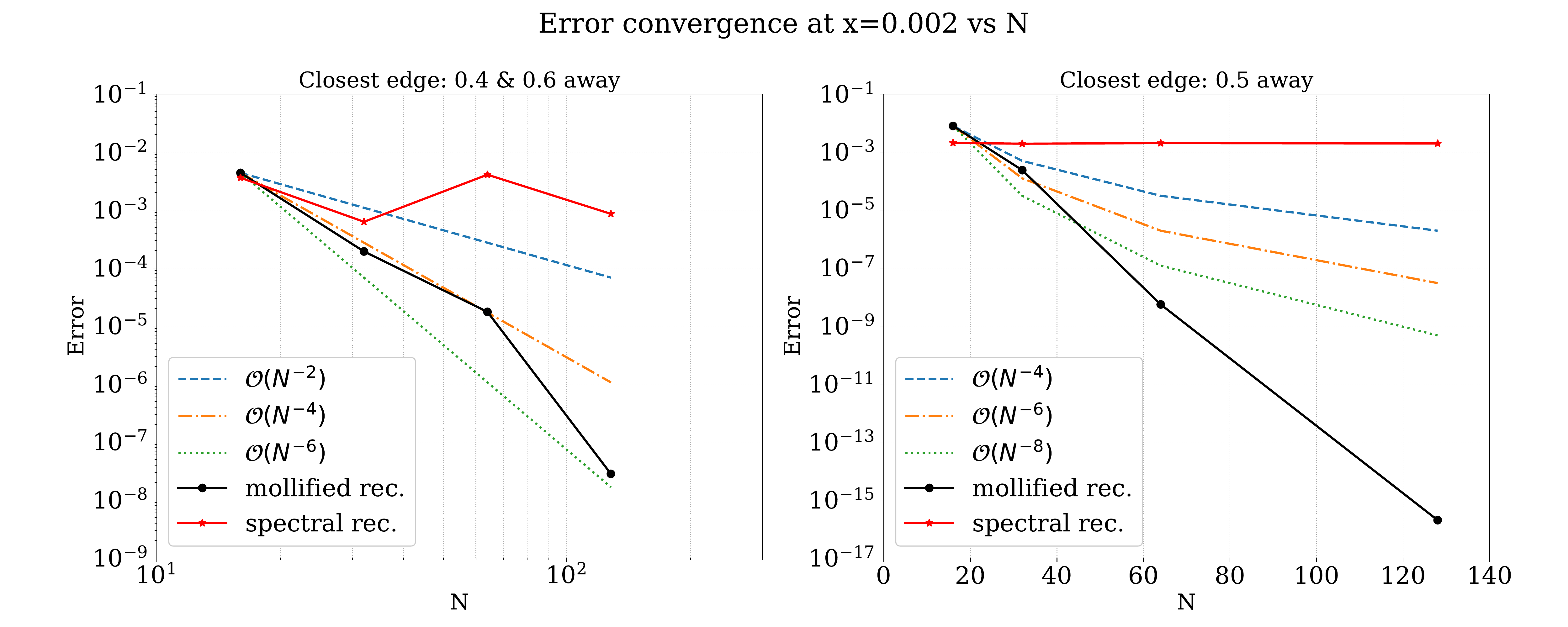}
  \caption{Difference between the analytic function and mollified
    spectral reconstruction at $x=0.002$ as a function of expansion
    order $N$. The point was chosen carefully, to avoid probing the
    spectral reconstruction at a collocation point. Cases under
    consideration are the same as in Figures
    \ref{fig:mollified:pointwise:error} and
    \ref{fig:mollified_unmollified}---two discontinuities within the
    domain (left) and a single discontinuity (right). Piecewise
    mollification seems to prove its exponential accuracy away from a
    single discontinuity (right), however suffers worse performance
    between the two of them (left). The spectral reconstruction
    behaves poorly in both cases, not exhibiting convergence at all.}
    %\TODO{The ordering of figures 13 and 14 is a bit weird given this
    %  figure. We might need to reverse the order. -JMM}}
  \label{fig:mollifier:convergence}
\end{figure*}

\begin{figure}
  \centering
  \includegraphics[width=0.45\textwidth]{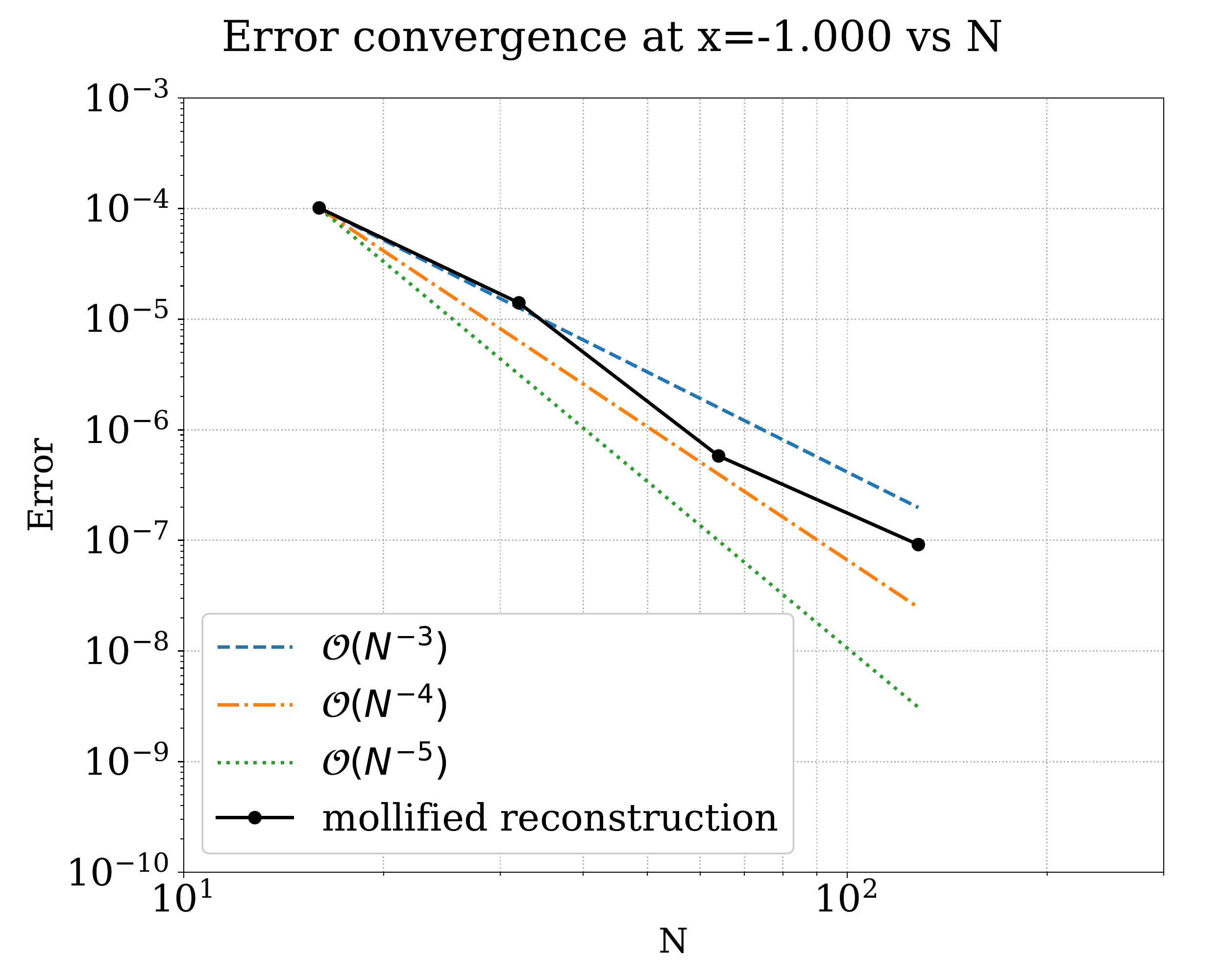}
  \caption{Difference between the analytic function and mollified
    spectral reconstruction at x=-1.00 domain boundary as a~function
    of expansion order $N$. The case under consideration is the same as
    in Figure \ref{fig:mollified:pointwise:error}---two
    discontinuities within the domain. Piecewise mollification seems
    to provide approximately 3rd order convergence at the domain
    boundary. The spectral reconstruction case is not included, as the
    domain boundary is one of the collocation points for
    Chebyshev-Gauss-Lobatto quadrature. By construction, this point is
    required to agree with the function exactly.}
  \label{fig:mollifier:convergence:boundary}
\end{figure}

Mollification, as discussed in Section
\ref{sec:intro:mollification:mollifiers} is an operation carried out
on the spectral reconstruction in real space. In our tests, we
approximate continuum mollification by discrete mollification on a
very fine, evenly-spaced grid of 500 points ($x$ values) within the
domain $[-1,1]$. For each $x$, we compute an appropriate adaptive
mollifier, based on the point position within the domain and the
location of the closest detected discontinuity. We then multiply the
Gibbs-polluted reconstruction by the mollifier at all $x$ and perform
discrete integration over the whole domain (e.g., via Chebyshev
integration or via a quadrature rule). This means we compute 500
integrals, which is potentially very expensive. For realistic
problems, a more efficient integration scheme is probably desirable.

\subsection{The Handling of Boundaries}
\label{sec:mollifiers:boundaries}

Non-periodic boundary conditions require careful handling of domain
boundaries. Lack of information outside of the domain influences
the integration result, with the effect most pronounced for grid
points closest to the edges of the domain at $x=-1.0$ and $x=1.0$. 
Every grid point outside of the domain effectively gives zero 
contribution to the discrete integral, creating an effect akin
to a smoothed discontinuity, discussed in Section
\ref{sec:intro:mollification:filters}.

In an attempt to improve convergence at the boundaries, we considered
different approaches such as the introduction of narrow buffer zones
outside of the domain and treating the boundaries as discontinuities
to renormalize the mollifiers. We first extended the information
beyond domain boundaries by adding a mirror reflection of the values
within $\delta(x)$ from the boundary \footnote{$\delta(x)$ is defined
  by \eqref{eq:def:delta} with $x=-1.0$ and $x=1.0$ for the respective
  boundaries.} to make sure we cover the width of essential support of
the mollifier and minimize the influence of missing values on the
discrete integral. We have chosen to mirror the values in order not to
introduce bias and only make use of the reconstructed function. This
approach yielded promising results, producing errors of order
$\mathcal{O}(10^{-3})$ for $N=32$ highest order in Chebyshev
expansion.

We then tried a different approach, using ideas derived from mollifier
modification across the discontinuity as described in Section 
\ref{sec:adhoc}. We required the mollifiers to reject information
beyond the domain by forcing their amplitude to zero across the boundary 
and renormalizing them to have unit mass. This solution resulted in
a better error behavior, with values of order $\mathcal{O}(10^{-5})$ for
$N=32$ highest order in Chebyshev expansion. It also had the advantage
of being entirely contained within the domain, therefore we decided
to implement it in our method.

\subsection{Results and error convergence}
\label{sec:mollifiers:results}

Piecewise mollification with the single-sided mollifiers illustrated
in Figures \ref{fig:mollifiers} and \ref{fig:mollifiers2D} has proved
effective in reducing Gibbs oscillations, while still preserving the
discontinuous character of the functions under investigation. Figure
\ref{fig:mollified32} shows the mollified result plotted against the
unprocessed spectral reconstruction and the original top hat function
computed for $N=32$ modes in the Chebyshev polynomial
expansion. Figure \ref{fig:mollified32:centre} presents the top hat
positioned well within the domain, while Figure \ref{fig:mollified32:edge}
illustrates a possible result of advection of the top hat waveform
such that majority of it has already managed to leave the domain.  In
both cases, the mollified reconstruction removes unwanted
oscillations, and retains a~(slightly smoothed out) discontinuity.

In order to quantify the method's performance we compute a pointwise
error across the domain by taking the difference between mollified
result and analytic function value at our 500 discrete $x$. We show
this error for the single discontinuity case in Figure
\ref{fig:mollified_unmollified} and the two discontinuity case in
Figure \ref{fig:mollified:pointwise:error}. The pointwise error
depends strongly on the position and number of
discontinuities. Nevertheless performance is quite satisfactory in all
cases we investigated. In the single discontinuity case, we also
compare the pointwise errors to the unmollified case. The improvement
via mollification is dramatic.

These relationships are further investigated in Figures
\ref{fig:mollifier:convergence} and
\ref{fig:mollifier:convergence:boundary}, where we plot the pointwise
difference between the mollified reconstruction and the analytic
function for $x=0.002$
%\footnote{point $x=0.0$ is one of the
%  collocation points for all $N$ in the figure, which are guaranteed
%  to coincide with the value of true underlying solution, as discussed
%  in \ref{sec:intro:projection}. However, mollification of the
%  spectral reconstruction does affect the result at all $x$ within the
%  domain, hence depriving $x=0.0$ of its special property. \TODO{Is
%    this footnote still relevant? -JMM}}
and
$x=-1.00$ respectively. We find that for a point away from a single
discontinuity the rate of error convergence exhibits exponential
behaviour, while it seems to be of order $\sim\mathcal{O}(N^{-4})$
between the two discontinuities. At the domain boundaries, regardless
of the solution character within it, the error seems to converge at
approximately 3rd order.

Our results are in qualitative agreement with
\revone{\cite{Tanner2006}}, where the author demonstrates exponential
error convergence away from discontinuities, nearly $\mathcal{O}(1)$
error at the discontinuity itself and similar magnitude of error
present at the domain boundaries.  However, it is important to
recognize that our numerical tests used different trial functions,
spectral projection and discrete integration methods, which do not
allow for a direct comparison between the two results.

\subsection{Stability of the method}
\label{sec:mollifiers:stability}

In some cases, the edge detection algorithm may fail to detect an edge
or locate it with significant error. (See for example Figure
\ref{fig:minmod:bad}.) In both cases the mollified solution still
retained its good behaviour away from the discontinuity, with the
increase in error influencing only the local environment of the
underlying real jump. In case of no edge detection close to the
boundary, the function would simply lose its discontinuous character
there and be smoothed out, as if we applied filtering to its spectral
coefficients. When the discontinuity was misplaced w.r.t to the
original jump, it would still remain discontinuous and without Gibbs
oscillations, however the $\mathcal{O}(1)$ error would be locally
extended in $x$. Further discussion of the robust character of the
method is presented in Section \ref{sec:hyperbolic:PDEs}.

\subsection{Computational Cost}
\label{sec:mollifiers:improvement}

Throughout our investigation, we performed the discrete integration at
every point using simple trapezoidal rule. To ensure that this
second-order error does not contribute to our analysis, we chose a
very fine grid for integration. Since one needs to perform one
integral per sampling point, this procedure introduces significant
computational cost---roughly the number of sampling points
squared. Reducing this high cost is an impediment to the successful
application of mollifiers to realistic problems and it requires
careful attention in the future.

\section{Hyperbolic PDEs}
\label{sec:hyperbolic:PDEs}

As a final test and as a proof-of-concept, we solve the
one-dimensional linear advection Equation
\begin{equation}
  \label{eq:def:advection:equation}
  \partial_t u - c \partial_x u = 0
\end{equation}
on the domain $x\in [-1,1]$ with periodic boundary conditions. We use
a discontinuous top-hat function like Equation \eqref{eq:def:tophat}
given as initial data. We solve Equation
\eqref{eq:def:advection:equation} using the method of lines and the
pseudospectral methods described in Section
\ref{sec:intro:projection}. We then post-process the solution via the
discontinuous mollifiers defined in Section \ref{sec:adhoc}. We detect
edges using Gelb's minmod method as described in Section
\ref{sec:edge:detection}. We recommend this same basic approach for
future applications of this method. By performing the pseudospectral
simulation, one can leverage the efficiency of spectral methods. Then,
when one wants to analyze the output, one can choose specific times
and snapshots to post-process and recover spectral accuracy via
mollification.

\subsection{Results for Discontinuous Mollification}
\label{sec:PDEs:mollification}

\begin{figure*}[t]
  \centering
    \large{Solution of the 1D advection equation with discontinuities \\
    post-processed with piecewise mollification}
  \par\smallskip
  \subfigure[]{
    \centering
    \label{fig:mollifier:near:center}
    \includegraphics[width=.43\textwidth]{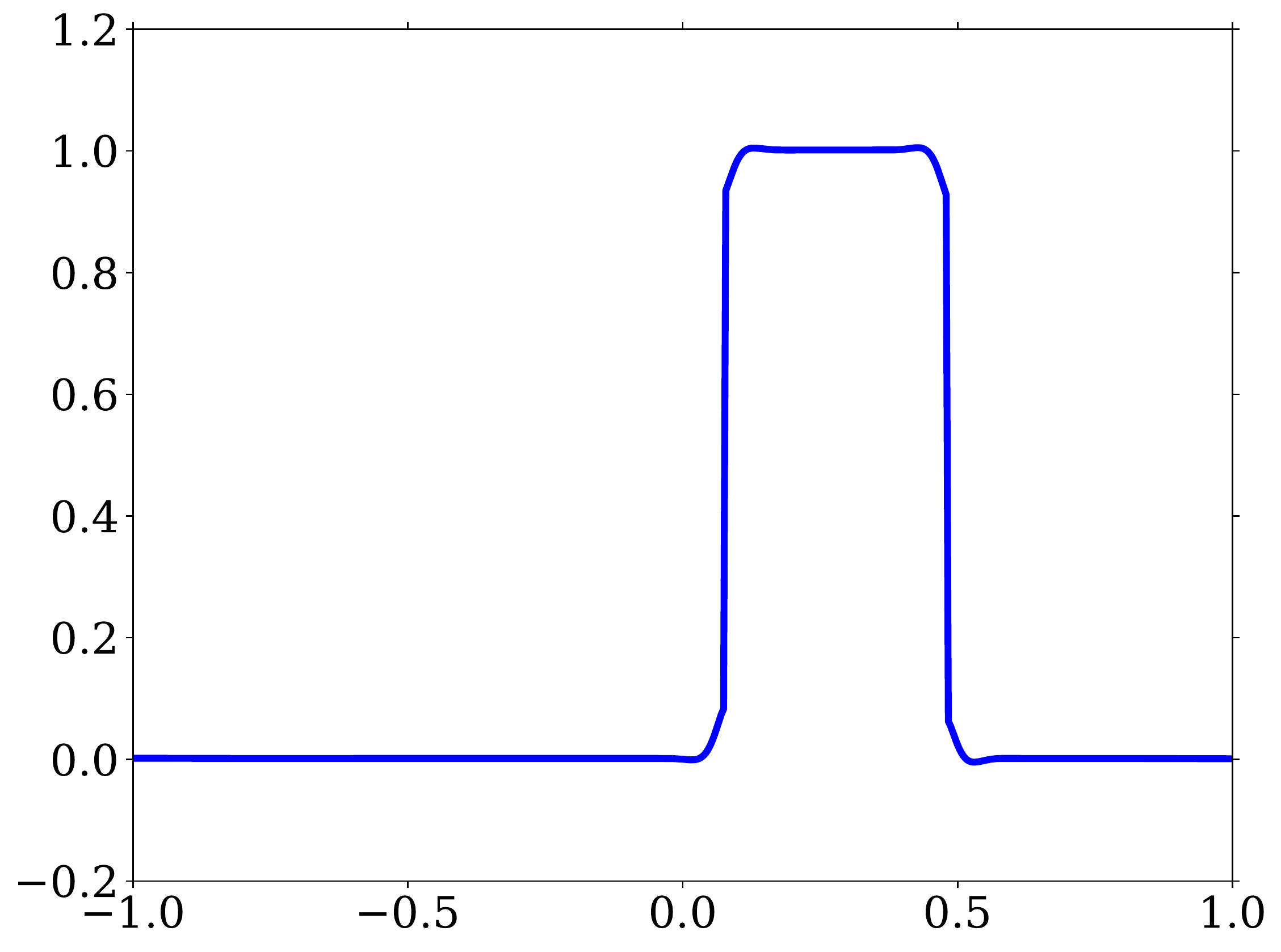}
  }
  \subfigure[]{
    \label{fig:mollifier:near:edge}
    \includegraphics[width=.43\textwidth]{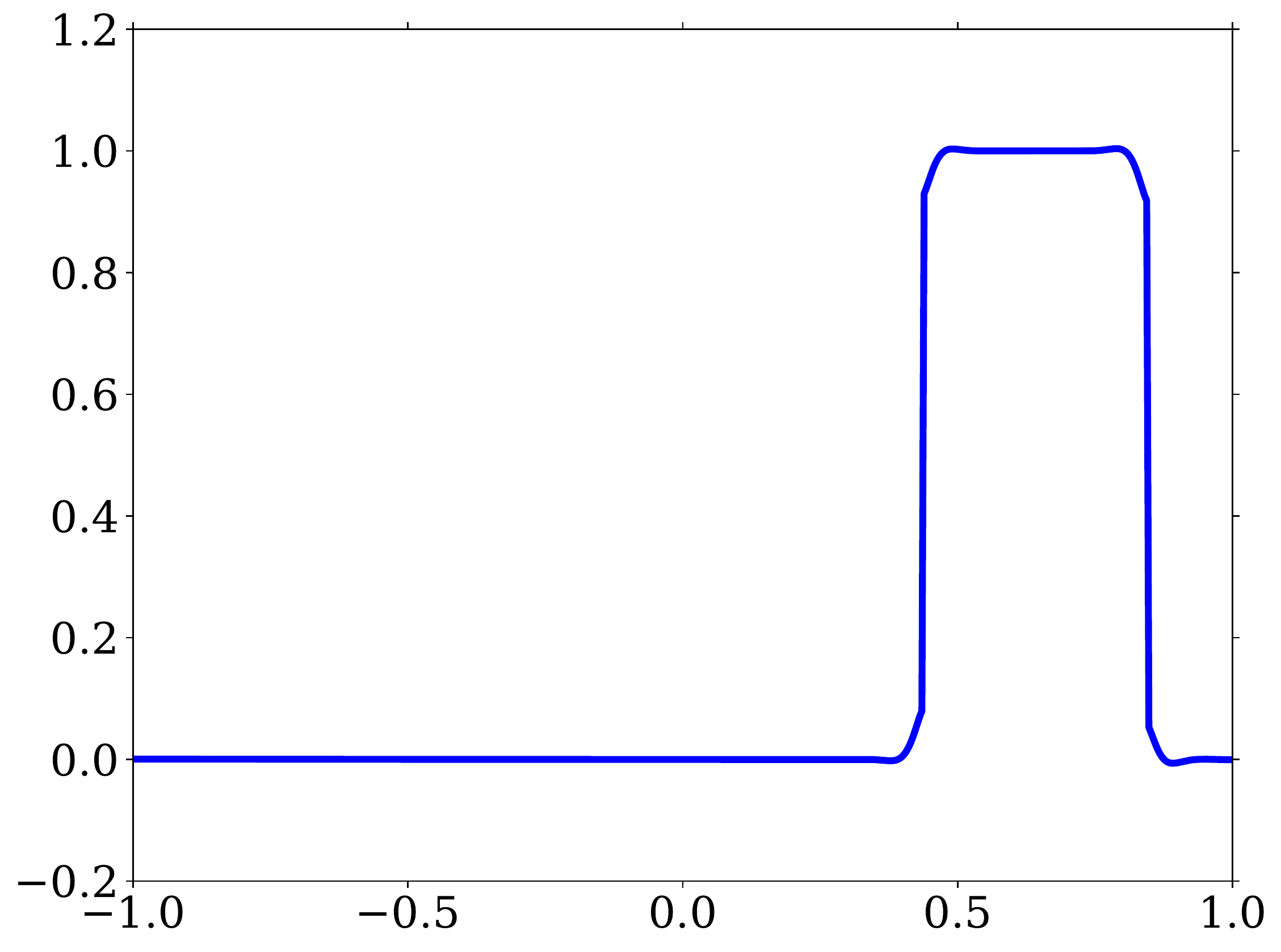}
  }
  \subfigure[]{
    \label{fig:mollifier:failing}
    \includegraphics[width=.43\textwidth]{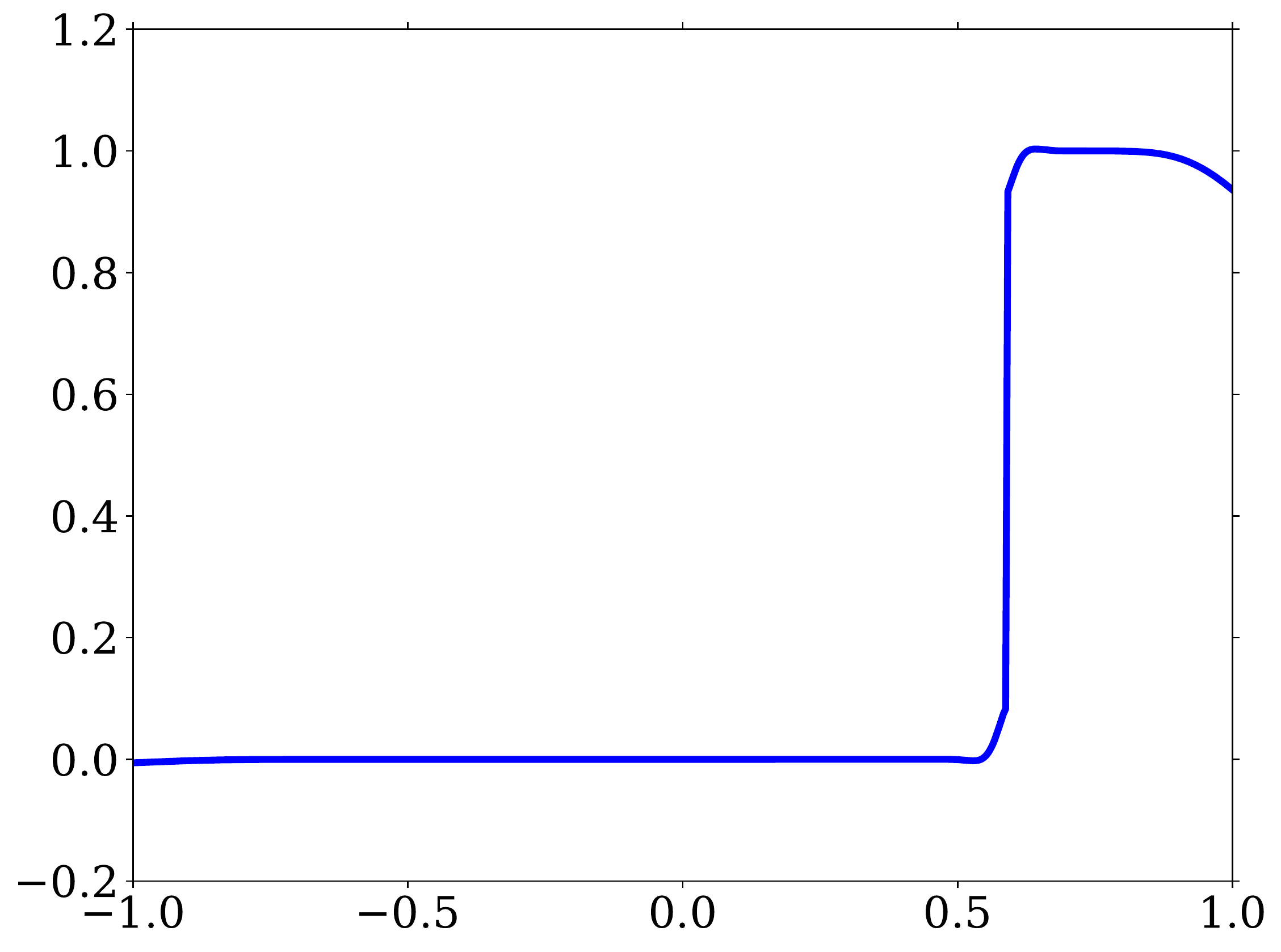}
  }
  \subfigure[]{
    \label{fig:mollifier:over:edge}
    \includegraphics[width=.43\textwidth]{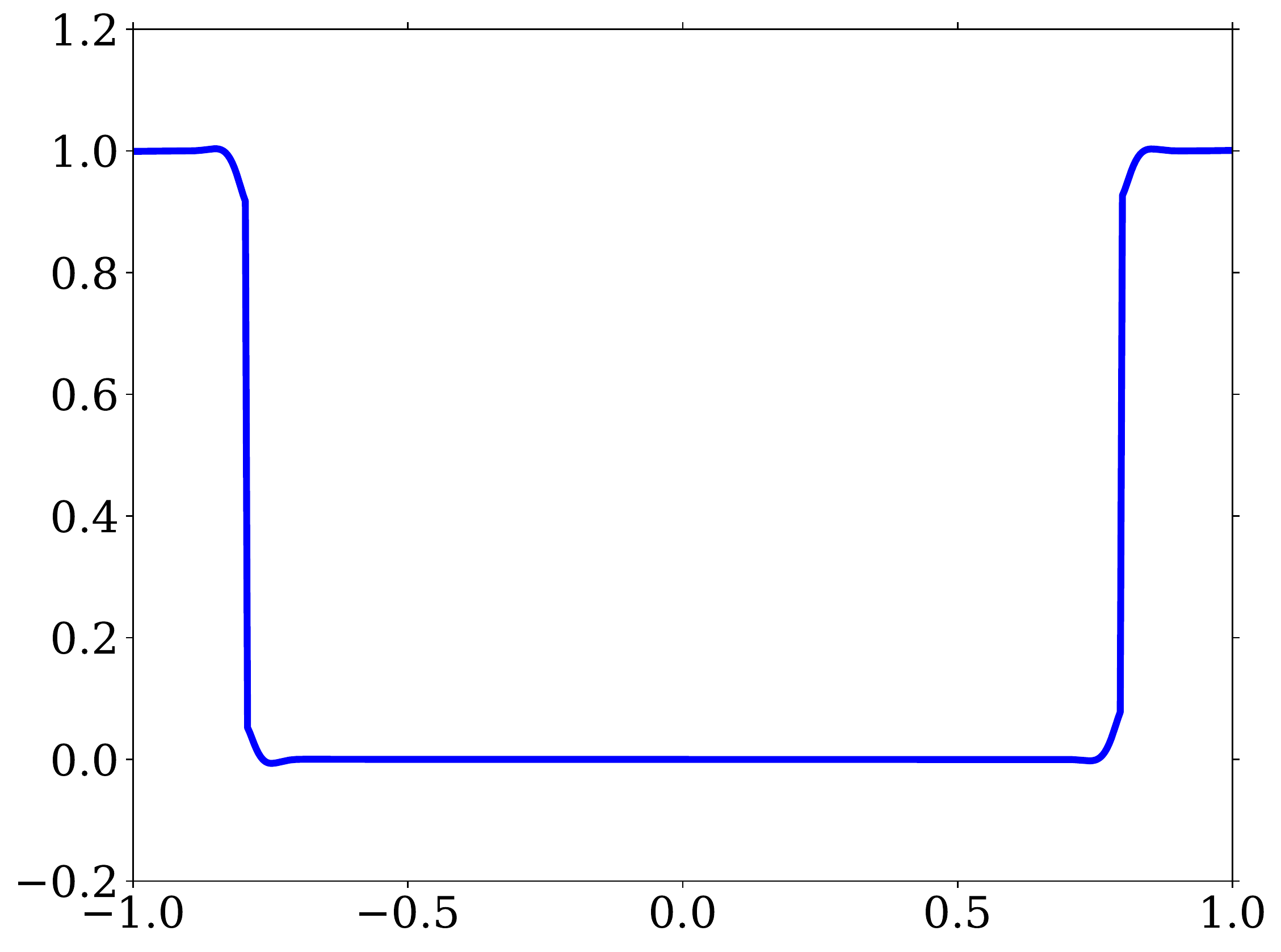}
  }
  \subfigure[]{
    \label{fig:mollifier:edge:at:bdry}
    \includegraphics[width=.43\textwidth]{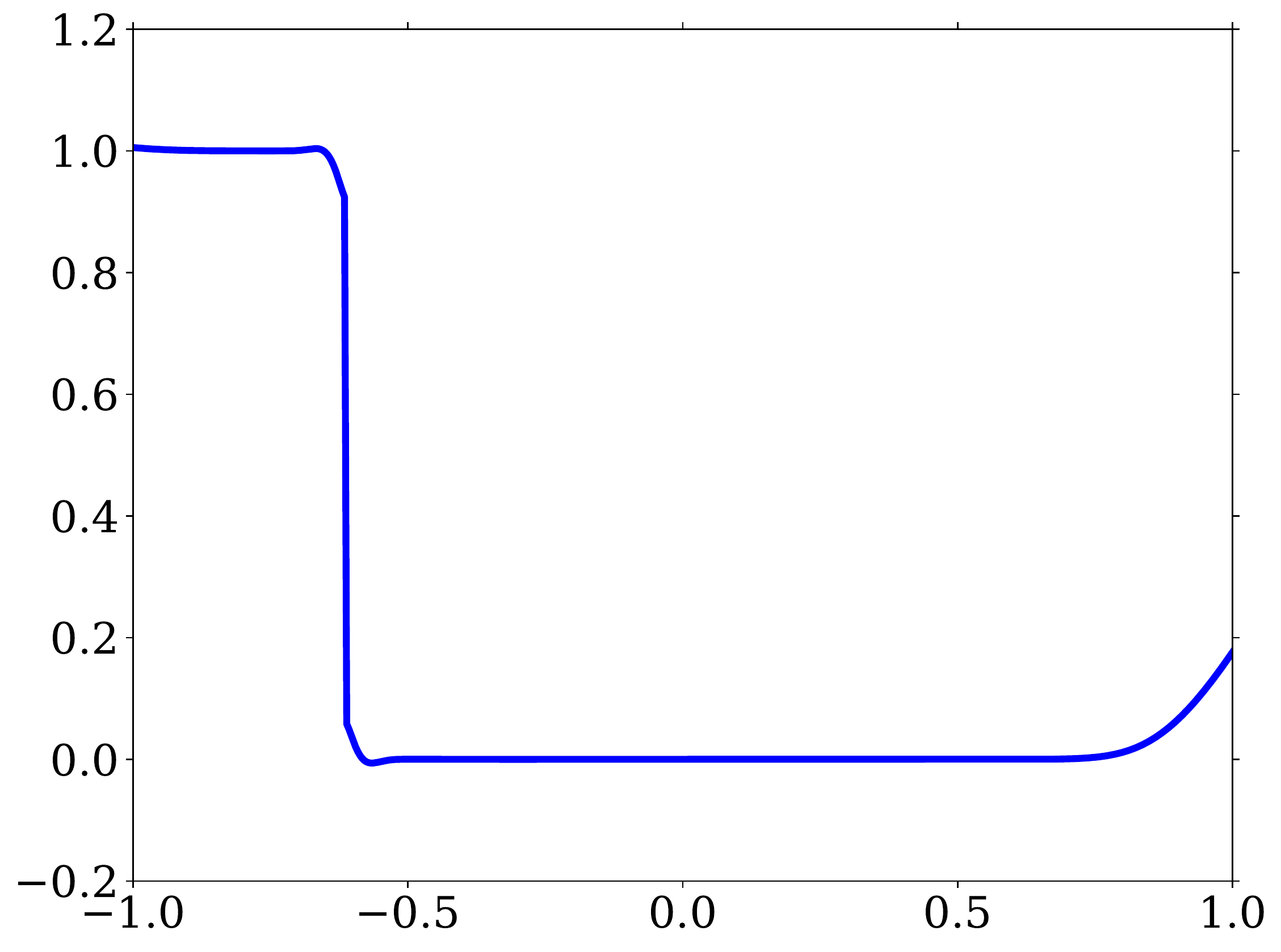}
  }
  \subfigure[]{
    \label{fig:mollifier:edge:calm}
    \includegraphics[width=.43\textwidth]{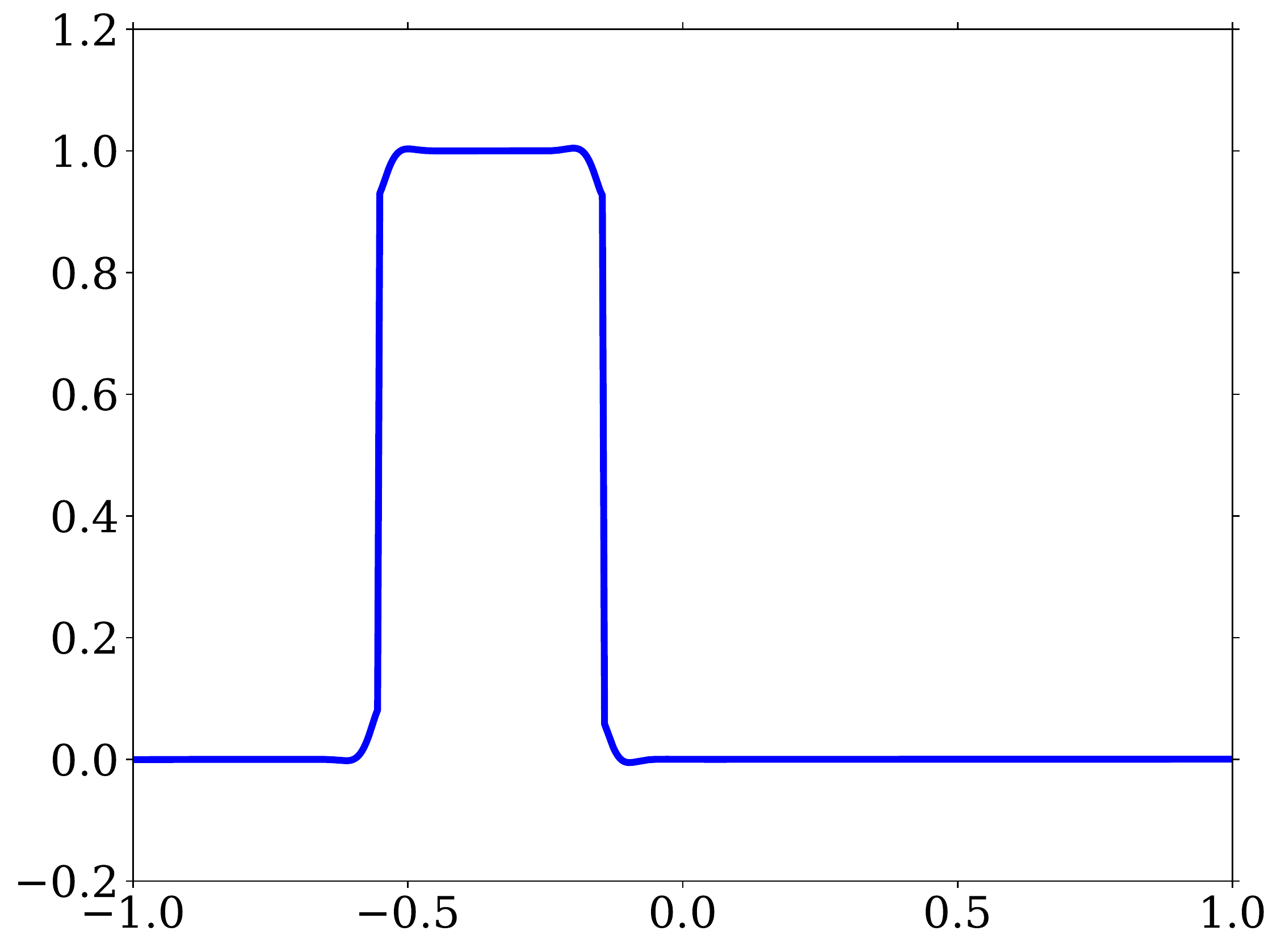}
  }
  \caption{Snapshots of the mollified reconstruction of an advected
    square wave on a periodic domain. The reconstruction suffers from
    significant errors when the wave reaches the domain boundary, when
    edge detection method fails. This, however does not entirely
    inhibit the recovery of the underlying solution. Read from left to
    right, top to bottom: \subref{fig:mollifier:near:center} the
    square wave begins well within the domain and the mollified
    reconstruction provides 4th order convergence presented in Figure
    \ref{fig:mollifier:convergence}; \subref{fig:mollifier:near:edge}
    as the wave approaches domain boundary the reconstruction is not
    affected; \subref{fig:mollifier:failing} the mollified
    reconstruction suffers significant error at the domain boundary,
    however does not prevent recovery of information between the
    detected discontinuity and the boundary;
    \subref{fig:mollifier:over:edge} the reconstruction remains stable
    as the wave crossed the periodic boundary;
    \subref{fig:mollifier:edge:at:bdry} the reconstruction is only
    affected in the direct vicinity of x=1.0 boundary;
    \subref{fig:mollifier:edge:calm} the reconstruction behaves just
    like in \subref{fig:mollifier:near:center}}
  \label{fig:advection:mollification}
\end{figure*}

We present snapshots of the mollified solution to Equation
\eqref{eq:def:advection:equation} in Figure
\ref{fig:advection:mollification}. When the discontinuities are far
from domain boundaries, the solution behaves as described in Section
\ref{sec:mollifiers:results}. There is a slight smearing of the
discontinuity (providing $\Ord{1}$ error), but otherwise the mollified
solution is extremely good. When a discontinuity is very near domain
boundaries, however, this smoothing increases and can spoil the
mollified solution near the discontinuity.\footnote{We emphasize that
  a failure of the mollifier \textit{does not} prevent the successful
  evolution of the system. All mollification is performed in
  post-processing.} Fortunately, this smearing is local, and the
global solution is still well-behaved.

% We believe this
% \textit{robustness} of the mollifier is critical for its successful
% application in the future. 

\subsection{On Robustness}
\label{sec:gegenbauer:PDEs}

For a numerical method to be applied successfully to real problems, it
must be \textit{robust}. The method must not only perform well in
ideal situations, but it must be capable of surviving non-ideal
cases, even if performance is reduced. In other words, a more
robust method fails less catastrophically than a less robust
method.

As demonstrated in the previous section, the convergence rate of our
one-sided mollifiers depends on the closeness to the discontinuity and
to the boundary of the domain. In the context of robustness, poor
convergence near the boundary of the domain is a failure. We argue
that mollification is \textit{robust} because this poor convergence is
a local effect and the solution is not spoiled in the rest of the
domain.

Here we argue that this robustness is a property of \revone{all} optimal
mollifiers, not just our one-sided ones.\footnote{It is likely a
  property of any local shock capturing technique} Moreover, we argue
that robustness cannot be assumed in general. To make this point, we
compare a naive implementation of \revone{optimal} mollifiers---as
described in section \ref{sec:intro:mollification:mollifiers}, without
our improvements in edge detection and without the one-sided property
introduced in section \ref{sec:adhoc}---to another shock-capturing
approach, the Gegenbauer reconstruction.

The Gegenbauer reconstruction, first proposed by Gottlieb and Shu
\cite{Gottlieb1992a} has received much attention as a post-processing
technique that enables spectral methods to resolve shocks and
discontinuities. The Gegenbauer polynomials $C^\alpha_n(x)$ are those
Jacobi polynomials which are orthogonal under the norm
\begin{equation}
  \label{eq:gegenbauer:norm}
  \braket{\phi,\psi} = \int \phi(x)\psi(x)(1-x^2)^{\alpha-1/2}dx.
\end{equation}
In this approach, smooth regions of a discontinuous spectral solution
are reprojected onto a finite-sized basis of a subset of the
Gegenbauer polynomials with fixed $\alpha$. $\alpha$~is chosen such
that the new basis is sufficiently ``different'' from the original
spectral basis. In this way, Gibbs oscillations can be removed up to
the discontinuity and a spectrally accurate, discontinuous solution
can be constructed. For more details on the Gegenbauer reconstruction,
see \cite{Gottlieb2011a} and references therein. For some examples of
``realistic'' one-dimensional applications of the Gegenbauer approach,
see
\cite{Gelb2000Enhanced,Gelb2001EnhancedSph,MeisterFilter}.\footnote{\cite{Gelb2001EnhancedSph}
  studies a two-dimensional system, but the Gegenbauer reconstruction
  is only applied to a one-dimensional test case.}

The Gegenbauer reconstruction thus fills the same role as our
discontinuous mollifiers. It is a post-processing technique that removes
the Gibbs oscillations from a~spectral solution of a PDE system and it
requires some method of finding the discontinuities, such as the one
described in Section \ref{sec:edge:detection}.

% It is therefore worth
% comparing the Gegenbauer  reconstruction to our mollifiers. 

% \footnote{Although the technique has some free
%   parameters in it, they can be chosen automatically. See, for
%   example, \cite{Gelb2005a,gelb2006robust}.}

\begin{figure*}[t]
  \centering
  \large{Robustness comparison between the Gegenbauer reconstruction \\
    and \revone{continuous mollification} in the solution of the 1D
    advection equation}
  \par\smallskip
  \subfigure[]{
    \centering
    \label{fig:gegenbauer:near:center}
    \includegraphics[width=.43\textwidth]{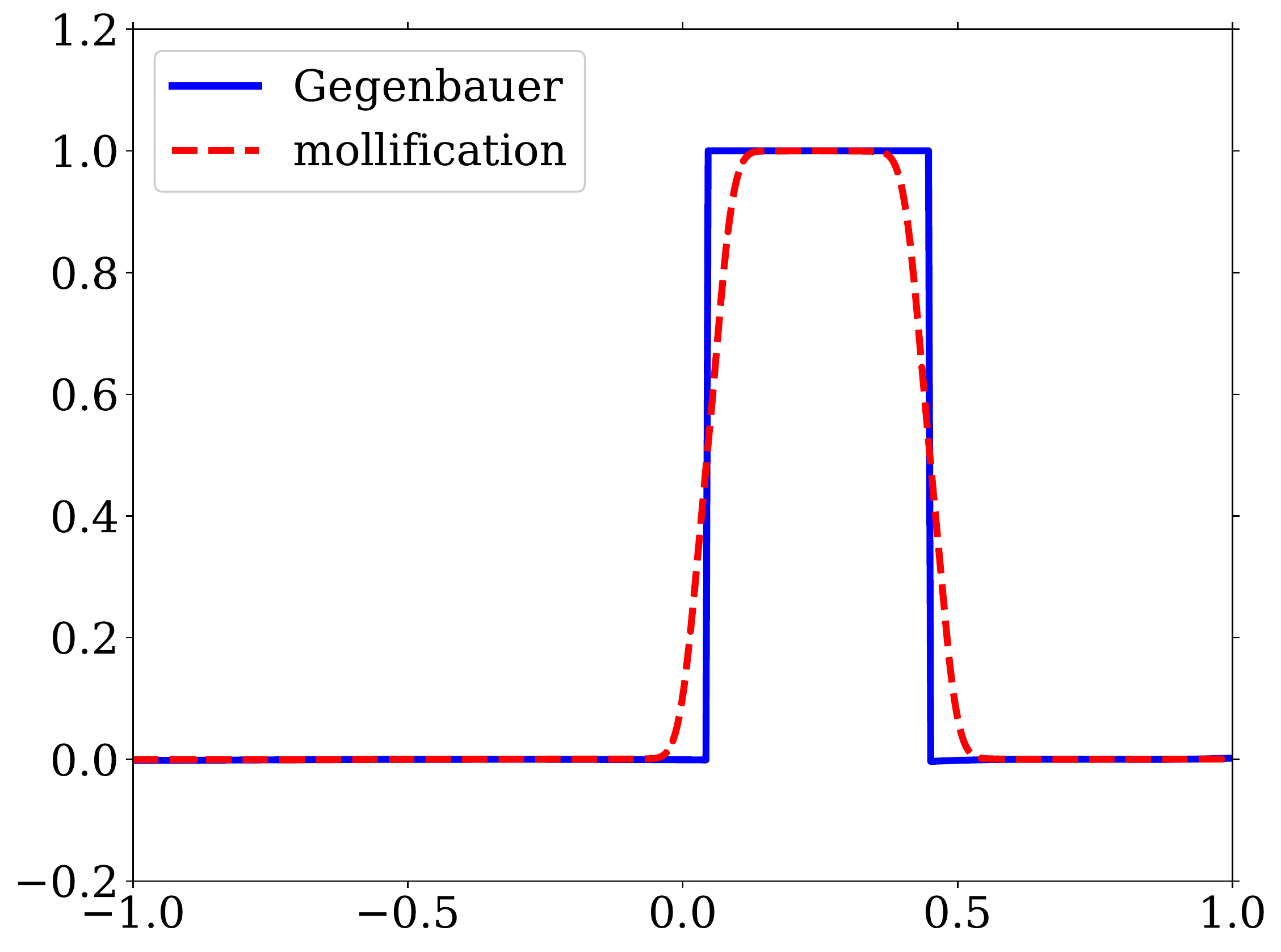}
  }
  \subfigure[]{
    \label{fig:gegenbauer:near:edge}
    \includegraphics[width=.43\textwidth]{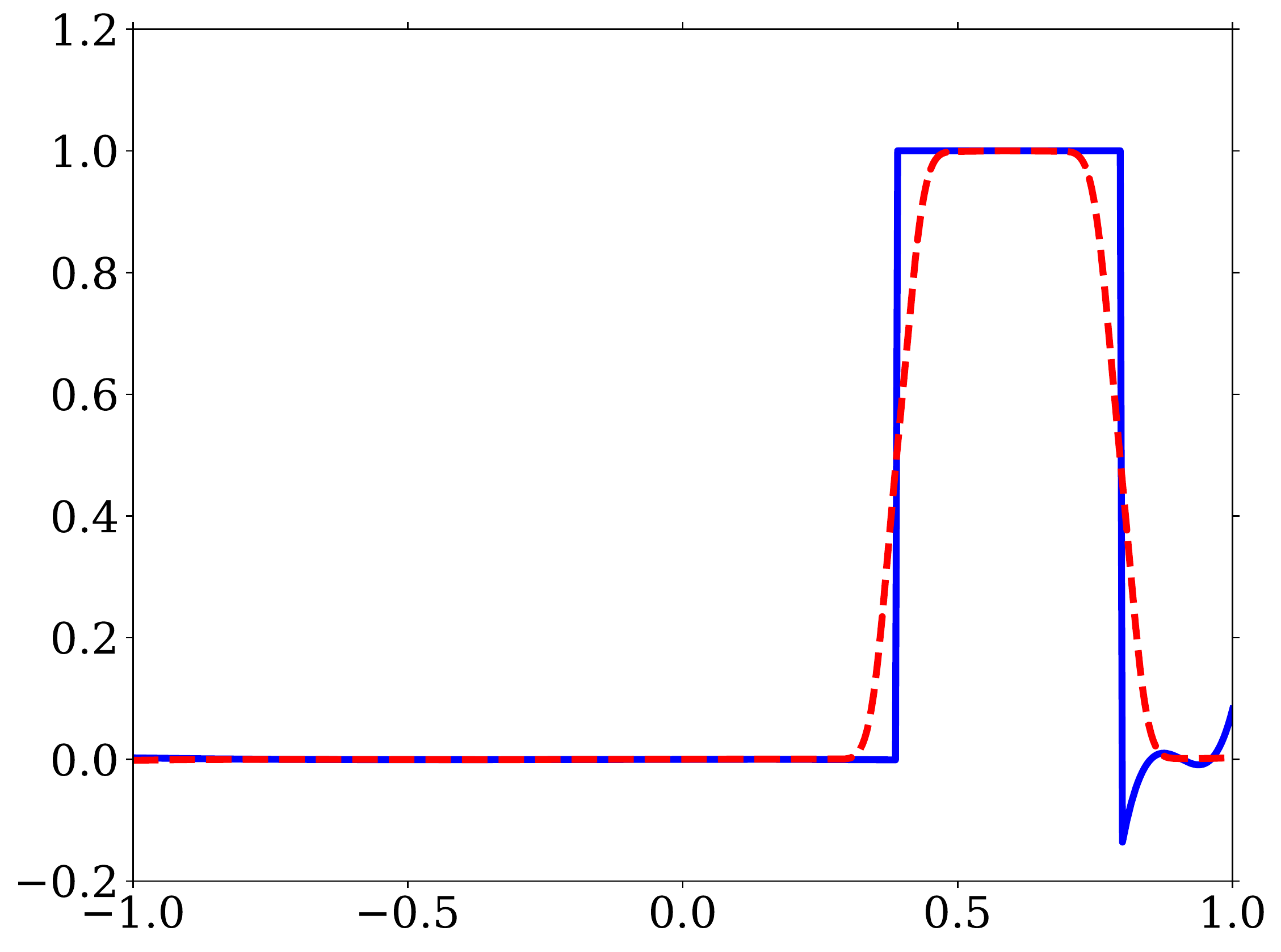}
  }
  \subfigure[]{
    \label{fig:gegenbauer:failing}
    \includegraphics[width=.43\textwidth]{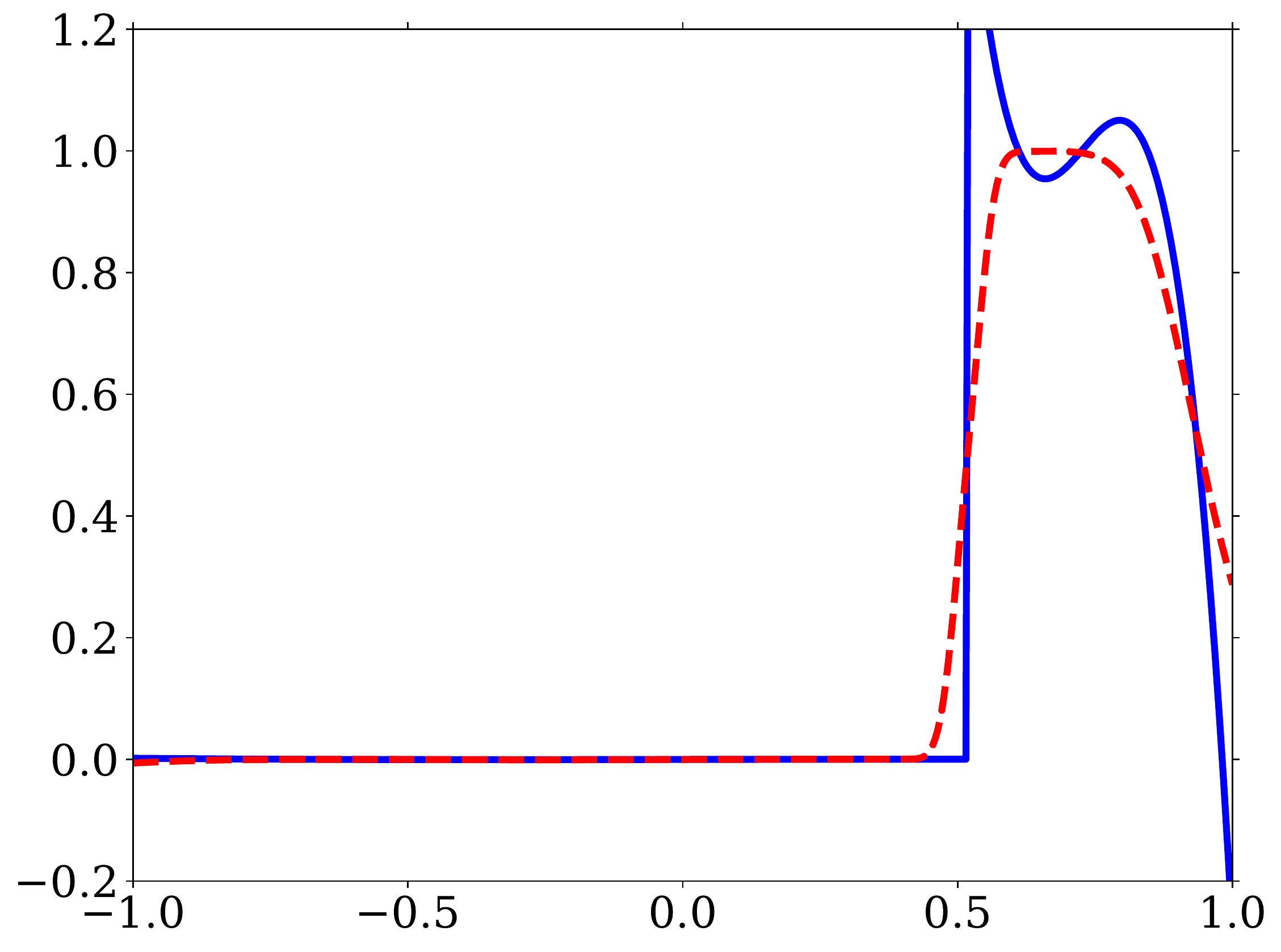}
  }
  \subfigure[]{
    \label{fig:gegenbauer:over:edge}
    \includegraphics[width=.43\textwidth]{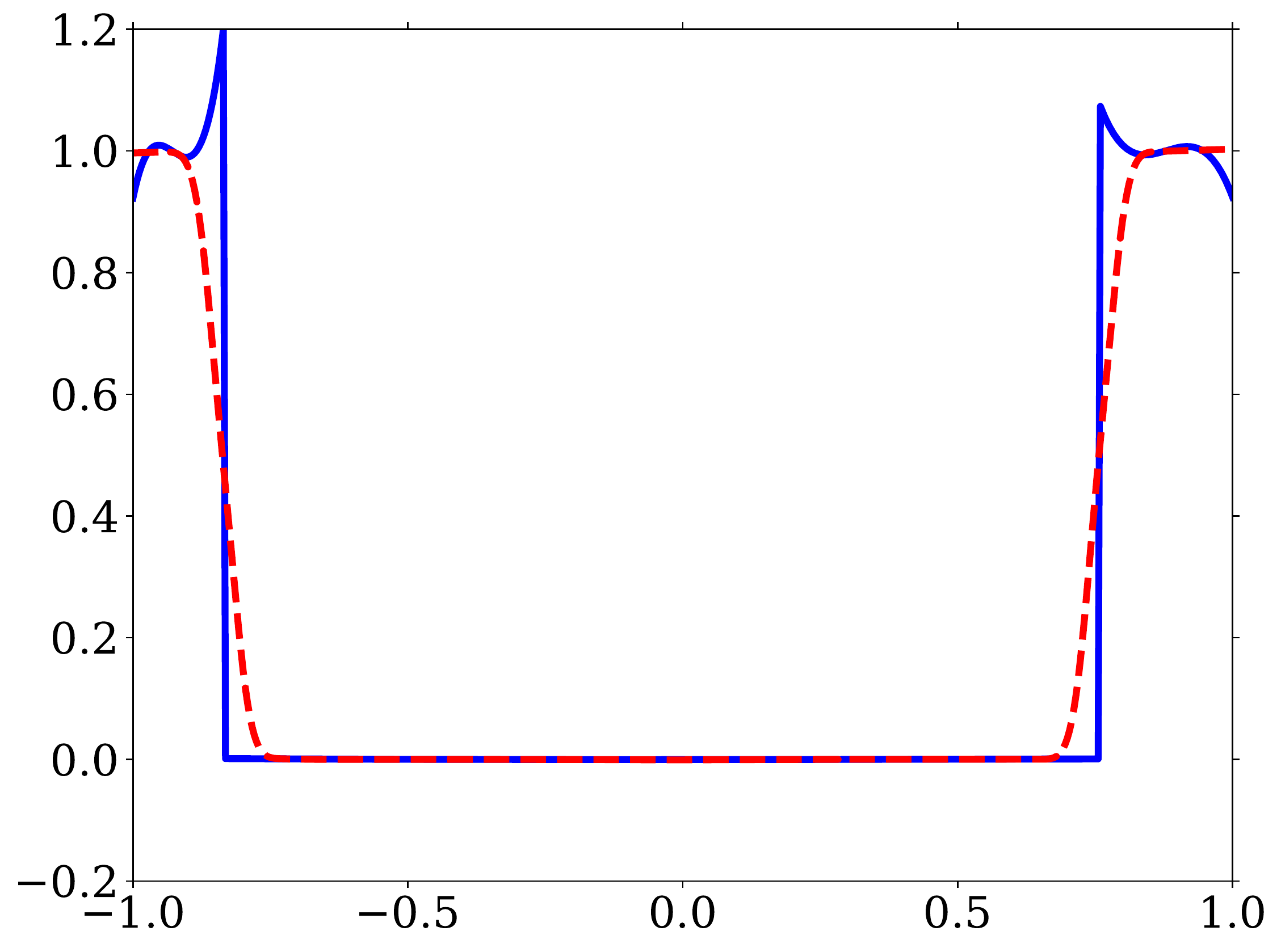}
  }
  \subfigure[]{
    \label{fig:gegenbauer:edge:at:bdry}
    \includegraphics[width=.43\textwidth]{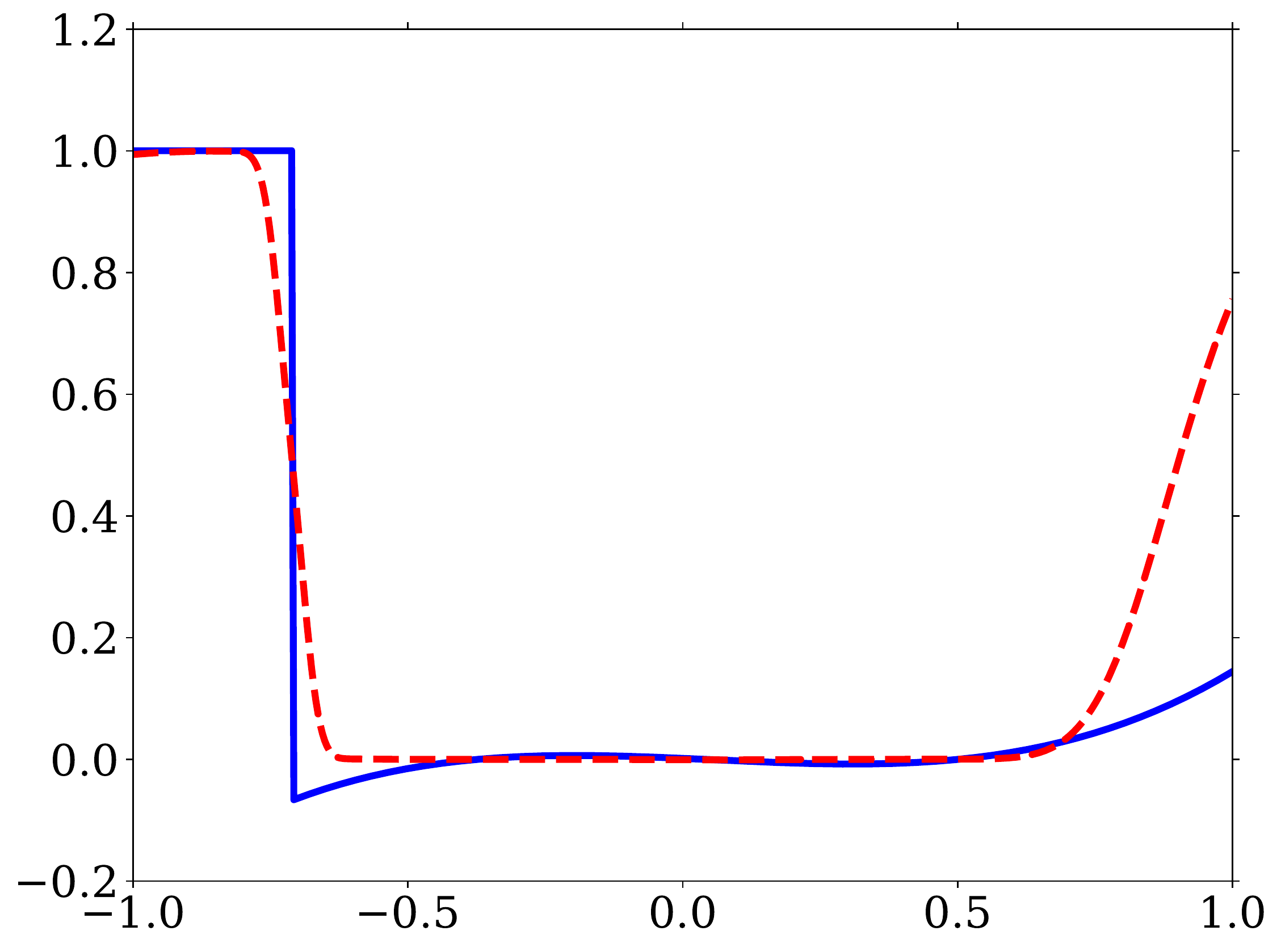}
  }
  \subfigure[]{
    \label{fig:gegenbauer:edge:calm}
    \includegraphics[width=.43\textwidth]{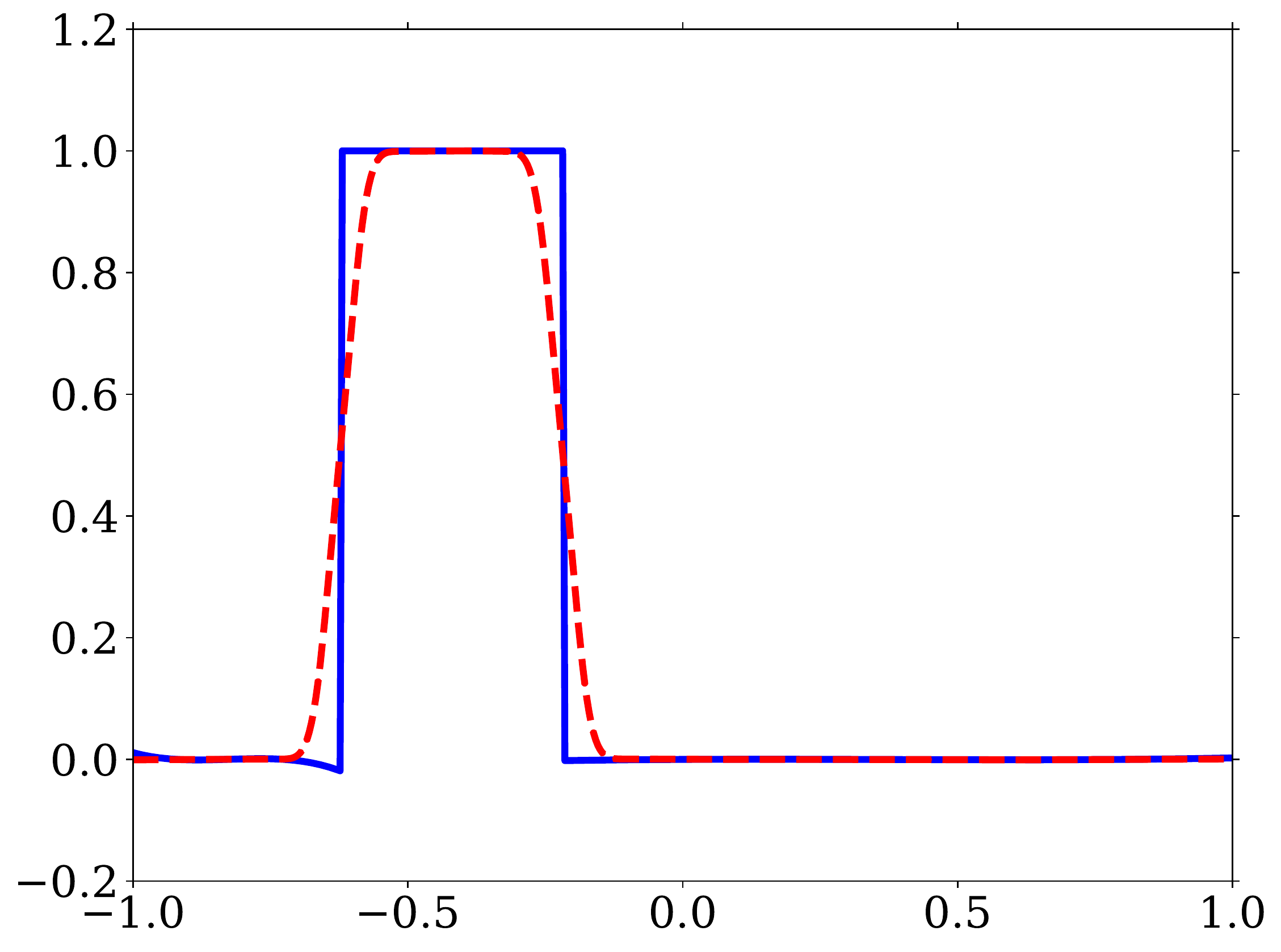}
  }
  \caption{Relevant snapshots of a Gegenbauer reconstruction (solid
    blue line) and a naive implementation of \revone{optimal
    mollifiers} (red dashed line) of an advected square wave on a
    periodic domain. The reconstruction fails when the wave reaches
    the periodic domain boundary. Read from left to right, top to
    bottom: \subref{fig:gegenbauer:near:center} the square wave begins
    near the center of the domain and the Gegenbauer reconstruction
    provides root-exponential convergence up to the point of
    discontinuity; \subref{fig:gegenbauer:near:edge} as the square
    wave approaches the edge, the reconstruction closest to the edge
    fails; \subref{fig:gegenbauer:failing} the Gegenbauer
    reconstruction fails entirely for the portion of the solution near
    the domain boundary; \subref{fig:gegenbauer:over:edge} the
    reconstruction stabilizes while the square wave crosses the domain
    boundary; \subref{fig:gegenbauer:edge:at:bdry} the reconstruction
    completely breaks down; \subref{fig:gegenbauer:edge:calm} the
    reconstruction recovers and converges exponentially far from the
    domain boundaries.}
  \label{fig:gegenbauer}
\end{figure*}

Figure \ref{fig:gegenbauer} compares a naive implementation of
mollifiers to the Gegenbauer reconstruction for the
advection equation \eqref{eq:def:advection:equation} using the
techniques described in Section \ref{sec:intro:projection} and the
method of lines. We use Equation \eqref{eq:def:tophat} as initial
data. Our Gegenbauer code is open source and the interested Reader can
find it in \cite{MillerGegenbauerCode}. The figure shows that the
reconstruction behaves well when the discontinuities are in the center
of the domain---indeed, the discontinuity is captured better than with
the mollifier. However, when the discontinuities are near domain
boundaries, the \textit{global} reconstructed solution breaks
down. The mollified solution does poorly near the domain boundaries,
but the global solution remains intact. In this way, the Gegenbauer
reconstruction is not \textit{robust}, and the consequences of this
failure of robustness are profound.

There are techniques designed to improve the reliability of the
Gegenbauer reconstruction---reducing the probability of
failure---which we did not explore. For example, see \cite{Gelb2005a},
\cite{gelb2006robust} and references therein. Our argument is not that
the Gegenbauer method is bad, or that it is worse than
mollifiers. Rather we argue that when the Gegenbauer reconstruction
fails, the failure is more damaging than when mollifiers fail.

In realistic settings, shocks and discontinuities may well be present
near domain boundaries. Moreover, in multiple dimensions, it may be
much more difficult to localize a discontinuity. We therefore believe
that \textit{robustness} of a method, i.e., how badly it damages the
solution when it fails, is a critical property that needs to be
prioritized when developing approaches to resolve the Gibbs
phenomenon.

\section{Comparison to Other Approaches}
\label{sec:other:approaches}

\revtwo{Mollifiers and the Gegenbauer reconstruction are by no means the only
way to deal with discontinuous functions. We would like to highlight a
few other approaches below, and discuss how they might compare to our
mollification scheme. Since we did not perform any numerical
experiments for these approaches, we cannot compare the robustness, as
described in section \ref{sec:gegenbauer:PDEs}, of these
methods. Nevertheless, we believe they merit mention. The techniques
discussed below have been extended to multiple dimensions.}

\revtwo{If the positions of the discontinuities are known a-priori, then a
discontinuous Galerkin (or any weak boundaries or penalty) method can
be very performant---one simply needs to place the physical
discontinuity at the element boundary
\cite{hesthaven2007nodal}. Although most discontinuous Galerkin
methods require filtering or artificial viscosity for stability, no
post-processing is needed and convergence is of order $\Ord{h^N}$,
where $h$ is the width of a discontinuous Galerkin element and $N$ is
the number of polynomials in the basis within an element. The
computational cost of a discontinuous Galerkin method is
$\Ord{h N^2}$. Since discontinuous Galerkin methods are highly
parallelizable and adaptive, this approach meets the needs of many
real science problems. For one example (of many) of this approach in
astrophysics, see \cite{FieldDGSelfForce}.}

\revtwo{The most naive application of discontinuous Galerkin methods and
similar approaches requires discontinuities whose positions are known
and unchanging. A more advanced approach might be Lagrangian (or
ALE), so that the mesh tracks the motion of discontinuities. This is a
challenging topic and an area of active research. For some examples of
these tricks, see \cite{NGUYEN20101,LIU201868}. Although it is very
promising, this approach carries with it significant infrastructure
and methods overhead, which we believe we avoid by focusing on
mollifiers.}

\revtwo{Another promising alternative post-processing
  technique is the moving least
squares method introduced by Lipman and Levin
\cite{LipmanLevinLeastSquares,AMIR201831}. In this approach, the error
introduced by the Gibbs phenomenon is modeled and fit via
least-squares. This model error can then be subtracted from the true
solution, recovering exponential convergence. A major advantage of the
moving least squares approach is that the cost scales linearly with
the number of collocation points, e.g., as $\Ord{N}$. This is because
the least-squares minimization occurs only locally around the
discontinuity and not over the whole domain. This is favorable
compared to mollification and the Gegenbauer reconstruction, which are
global and scale as $\Ord{N^2}$ in cost. (This is of course, only for
the post-processing step. The evaluation of the time-evolution
operator of a spectral method is $\Ord{N^2}$, regardless.) We could
not find any examples of this approach being used in practice,
however, we believe it merits further study.}

% \TODO{REREAD ME AND FIX ME}

\section{Concluding Thoughts}
\label{sec:conclusion}

% \TODO{Clean up the conclusion. -JMM}

% \TODO{DISCUSS EXTENSION TO MULTI-D}

In the course of this work, we have investigated the efficacy of
several techniques that detect discontinuities in spectral data and
that help reduce the error introduced by the associated Gibbs
phenomenon.
We have found that, although these techniques are very promising, they
require careful implementation for practical applications. The
spectral edge detection developed by Gelb is a good example of
this. Theoretical accuracy is attained in ideal situations, such as
for a discontinuity in the center of the domain, but spurious
oscillations spoil the solution when the discontinuity is near a
domain boundary. (See, e.g., figures \ref{fig:edgeloc:centre} and
\ref{fig:edgeloc:boundary}.)

Moreover, the best technique in theory may be inferior in practice. In
ideal cases, the Gegenbauer reconstruction provides uniform
root-exponential convergence (in the number of modes of the original
system) to a true solution, up to and including discontinuities
\cite{Gottlieb1992a}. In our own experiments, we found the Gegenbauer
reconstruction to perform excellently in ideal situations, but that
when it failed---for example, when the discontinuity was poorly
localized (as is the case with Gelb's edge detection near domain
boundaries)---it failed catastrophically. \revone{Optimal} mollifiers also
failed in these situations, but the failure did not spoil the global
solution. (See, e.g., figures \ref{fig:gegenbauer} and
\ref{fig:advection:mollification}.)

% We believe our experience aligns
% well with Boyd's argument.

% However, Boyd demonstrated that as the number of
% Gegenbauer modes increases, convergence can easily be destroyed
% \cite{BoydTroubleWithGegenbauer}.

As a result of our experiments, we have improved upon both the edge
detection developed by Gelb and collaborators
\cite{Gelb1999a,Gelb2001a,Gelb2008a} and the mollifiers developed by
Tadmor, \revone{Tanner}, and collaborators
\cite{gottlieb1985recovering}, \cite{Tadmor2002},
\cite{Tanner2006}. In particular, we have developed several techniques
for handling erroneous behaviour at the boundaries of the domain and,
more importantly, we have introduced mollifiers that vanish outside
the region of smoothness. These \textit{one-sided} mollifiers allow
for the recovery of truly discontinuous solutions in a way that is
resistant to perturbation. We believe our improved versions of these
techniques are a promising path towards using high-order spectral
methods in production simulations of non-smooth problems. However,
there are many obstacles that must be overcome. The study presented
here is preliminary---it is only in 1D and covers relatively simple,
linear systems.

In higher dimensions, the increased dimensionality of the
discontinuity itself poses a problem. In one dimension, a
discontinuity is localized to a point. However, in two and three
dimensions it is a line and a surface respectively. Worse,
the extra dimensionality allows discontinuities to have complex
geometric and topological structure.\footnote{And of course, in many
  dimensions, as in one, there may be multiple disjoint discontinuities.} The
obvious extension of our one-dimensional results to this setting is
via a Cartesian product. However, it is not clear how effective this
approach will be.

Nonlinear systems provide another challenge. In the linear case, the
error introduced by the Gibbs phenomenon simply advects across the
grid. However, in nonlinear systems, aliasing error will translate
this error into ``physical'' source terms and can drive an
instability. Tadmor's spectral viscosity method provides stability and
convergence criteria for spectral solutions of these nonlinear systems
\cite{Tadmor1990}. However, it is not clear how effective this
technique will be when multiple kinds of physics and solution methods,
for example gravity and radiation transport, are added to the
calculation.

Another issue with nonlinear systems is that they have regimes of
validity. For example, the Gibbs phenomenon may drive fluid density or
pressure to a non-positive value. One may need to play tricks (such as
an artificial atmosphere) to preserve stability. And these tricks may
damage global convergence.

On the other hand, nonlinear systems also provide physical measures of
success and physical tools which can be used to detect shocks and
discontinuities, perhaps augmenting the edge detection discussed
here. One could, for example, look at generalized Riemann conditions
or entropy production for contact discontinuities and shock
discontinuities respectively.

Finally, mollification efficiency is a challenge. Mollifying a
spectral reconstruction requires one integral per sample point. If the
quadrature points of numerical integration are the sample points, this
translates into a~computational cost quadratic in the number of sample
points. This cost is somewhat mitigated by the fact that mollification
is performed as a post-processing step and does not need to be applied
at every time step. Nevertheless, an efficient and adaptive
integration scheme is required.\footnote{We note that the fact that
  the mollifier depends on space, and in fact changes dramatically
  when one moves across a discontinuity, precludes its evaluation in
  the spectral domain. Therefore, we can not take advantage of more
  efficient FFT methods to improve efficiency.}

We do not believe these difficulties are insurmountable. However, they
require careful study. Indeed, we believe they are promising avenues
of future research. The extreme efficiency of spectral methods mean
that, if these difficulties are resolved, calculations that are
currently performed on a supercomputer could be performed on a
desktop. These potential gains make this extra care worthwhile.

\section{Acknowledgments}
\label{sec:acknowledgements}

We thank Roland Haas and Oleg Korobkin for their insightful comments
on the draft. JMM thanks David Radice and our anonymous peer reviewer
for guiding him to the appropriate literature.

The authors acknowledge support from the Natural Sciences and
Engineering Research Council of Canada (NSERC).
The research was also supported by the Perimeter Institute for
Theoretical Physics. Research at Perimeter Institute is supported by
the Government of Canada through Industry Canada and by the Province
of Ontario through the Ministry of Research and Innovation.

This work was supported by the US Department of Energy through the
Laboratory Directed Research and Development Program of Los Alamos
National Laboratory at Center for Nonlinear Studies under project
number 20170508DR. Los Alamos National Laboratory is operated by Triad
National Security, LLC, for the National Nuclear Security
Administration of U.S. Department of Energy (Contract
No. 89233218CNA000001).

JMM acknowledges the U.S. Department of Energy Office of Science and
the Office of Advanced Scientific Computing Research via the
Scientific Discovery through Advanced Computing (SciDAC4) program and
Grant DE-SC0018297.

We are grateful to the countless developers contributing to open
source projects on which we relied in this work, including Python
\cite{rossumPythonWhitePaper}, numpy and scipy \cite{numpy,scipyLib},
and Matplotlib \cite{hunterMatplotlib}.

%\TODO{Any other open source projects you want to thank?}

% \TODO{HPC systems+allocations?}

\appendix

% \section{Some Appendix}
% \label{sec:some:appendix}
% 
% this is some appendix.
% 
% \begin{widetext}
%   \section{Wide Appendix}
%   \label{sec:wide:appendix}
%   this is how we make things
%   SUUUUUUUUUUUUUUUUUUUUUUUUUUUUUUUUUUUUUUUUUUUUUUUUPER wide so that
%   they extend beyond a single column.
% \end{widetext}

%bibliography
\bibliography{spectral-methods-in-presence-of-discontinuities}
\bibliographystyle{habbrv}

\end{document}